\RequirePackage{fix-cm}
\documentclass[smallextended]{svjour3}       % onecolumn (second format)
\smartqed  % flush right qed marks, e.g., at end of proof

\usepackage{amsmath,amsthm,amssymb}
\usepackage{mathtools}
\usepackage{csquotes}
\usepackage{mathtools,calc}
\usepackage{pgfplots}
\usepackage{xfrac}
\usepackage{framed}
\usepackage[ruled,vlined]{algorithm2e}
\usepackage{adjustbox}
\usepackage{makecell}
\usepackage{multirow}
\usepackage{bbding}
\usepackage{tikz}
\DeclareMathAlphabet{\mathbcal}{OMS}{cmsy}{b}{n}
\usepackage{appendix}

\newcommand{\diff}{\text{d}}           
\newcommand{\flux}{\mathcal{F}}

\newcommand{\numU}{\mathcal{U}}
\newcommand{\numS}{\mathcal{S}}

\usepackage[T2A,T1]{fontenc}
\newcommand{\Zhe}{\mbox{\usefont{T2A}{\rmdefault}{m}{n}\CYRZH}}
\usepackage{graphicx}
\usepackage{subcaption}
\usepackage{float}
\usepackage[framemethod=tikz]{mdframed}

\usepackage{url}
\usepackage{xcolor}
\definecolor{newcolor}{rgb}{.8,.349,.1}

\usepackage{lineno} 

\begin{document}

\title{A modified equation analysis for immersed boundary methods based on volume penalization: applications to linear advection-diffusion and high-order discontinuous Galerkin schemes %\thanks{Grants or other notes %about the article that should go on the front page should be %placed here. General acknowledgments should be placed at the end of the article.} 
}
\titlerunning{Modified equation analysis for IBM-VP DGSEM}

\author{Victor J. Llorente \and Jiaqing Kou \and Eusebio Valero \and Esteban Ferrer }


\institute{V. J. Llorente \at
              ETSIAE-UPM-School of Aeronautics, Universidad Politécnica de Madrid, Plaza Cardenal Cisneros 3, E-28040 Madrid, Spain\\ \email{victorjavier.llorente@upm.es}
           \and
            J. Kou (corresponding author) \at ETSIAE-UPM-School of Aeronautics, Universidad Politécnica de Madrid, Plaza Cardenal Cisneros 3, E-28040 Madrid, Spain\\
              \email{jiaqingkou@gmail.com}           %  \\
           \and
           E. Valero \at
              ETSIAE-UPM-School of Aeronautics, Universidad Politécnica de Madrid, Plaza Cardenal Cisneros 3, E-28040 Madrid, Spain\\
              Center for Computational Simulation, Universidad Politécnica de Madrid, Campus de Montegancedo, Boadilla del Monte, 28660 Madrid, Spain\\
           \and
           E. Ferrer \at
              ETSIAE-UPM-School of Aeronautics, Universidad Politécnica de Madrid, Plaza Cardenal Cisneros 3, E-28040 Madrid, Spain\\
              Center for Computational Simulation, Universidad Politécnica de Madrid, Campus de Montegancedo, Boadilla del Monte, 28660 Madrid, Spain\\
}

\date{Received: date / Accepted: date}

\maketitle

\begin{abstract}
The Immersed Boundary Method (IBM) is a popular numerical approach to impose boundary conditions without relying on body-fitted grids, thus reducing the costly effort of mesh generation. To obtain enhanced accuracy, IBM can be combined with high-order methods (e.g., discontinuous Galerkin). For this combination to be effective, an analysis of the numerical errors is essential. In this work, we apply, for the first time, a modified equation analysis to the combination of IBM (based on volume penalization) and high-order methods (based on nodal discontinuous Galerkin methods) to analyze \textit{a priori} numerical errors and obtain practical guidelines on the selection of IBM parameters. The analysis is performed on a linear advection-diffusion equation with Dirichlet boundary conditions. Three ways to penalize the immerse boundary are considered, the first penalizes the solution inside the IBM region (classic approach), whilst the second and third penalize the first and second derivatives of the solution. We find optimal combinations of the penalization parameters, including the first and second penalizing derivatives, resulting in minimum errors. We validate the theoretical analysis with numerical experiments for one- and two-dimensional advection-diffusion equations. 

\end{abstract}

%\tableofcontents

%\linenumbers

\section{Introduction}\label{sec:Intro}
Despite successful applications of Immersed Boundary Method (IBM) in simulating complex flows \cite{mittal2005immersed,huang2019recent} and fluid-structure interaction problems \cite{sotiropoulos2014immersed,kim2019immersed,griffith2020immersed}, understanding and controlling numerical errors in the IBM approach remains a challenge. IBM refers to a group of numerical strategies that handle the boundary condition when the solid is immersed in the computational domain, avoiding body-fitted meshes and enabling the use of simple meshes (e.g., Cartesian or Octree). The IBM approach originates from the idea of Peskin \cite{peskin1972IBM}, where singular forces represented by delta functions were positioned at solid boundaries to mimic the effect of physical boundaries. Since IBM reduces the complexity of mesh generation and handles moving boundaries efficiently, it has received a lot of attention over the past few decades. In general, IBM treatment can be achieved using the cut cell approach \cite{ye1999accurate,udaykumar2001sharp,fidkowski2007triangular,sticko2019higher} or by introducing source terms, e.g., ghost cell \cite{majumdar2001rans}, projection method \cite{taira2007immersed}, direct forcing \cite{fadlun2000combined,luo2012numerical} or volume penalization \cite{angot1999penalization,kolomenskiy2009fourier}. Although the cut-cell approach shows better convergence properties, the extension to moving boundaries and the treatment of different types of cut-cells remain challenging. A more flexible approach is the IBM based on Volume Penalization (VP). The latter shows advantages in robustness, ease of implementation, and theoretical convergence estimates~\cite{arquis1984conditions,angot1999penalization}.

Volume penalization is a classic IBM treatment based on modeling the solid as a porous medium with low permeability \cite{kadoch2012volume,sakurai2019volume}. The method imposes the boundary condition by introducing a source term (or penalty term) to the computational nodes located inside the solid. This approach dates back to the work of Courant~\cite{courant1943variational}, where a penalty method was used to transform constrained optimization problems into a constraint-free problem. Volume penalization methods for the Navier-Stokes equations were first proposed by Arquis and Caltagirone~\cite{arquis1984conditions} with a Brinkman-type penalization for the momentum equations. After that, Angot et al.~\cite{angot1999penalization} and Carbou and Fabrie~\cite{carbou2003boundary} proved the convergence of volume penalization, showing that as the penalization parameter $\eta$ approaches $0$, the model error converges if non-slip boundary conditions are considered. Subsequently, the volume penalization was extended to allow Neumann boundary conditions \cite{kadoch2012volume,sakurai2019volume} and Robin boundary conditions \cite{ramiere2007general}, as well as spatially varying Neumann and Robin boundary conditions \cite{thirumalaisamy2021handling}. The volume penalization method was extended to compressible flows by Liu and Vasilyev~\cite{liu2007brinkman}, Brown-Dymkoski et al.~\cite{brown2014CBVP}, Abgrall et al.~\cite{abgrall2014IBM} and Abalakin et al.~\cite{abalakin2016immersed}. This method has been applied to a variety of problems, including flapping wings~\cite{kolomenskiy2009fourier}, two-phase flows~\cite{horgue2014penalization}, aeroacoustics \cite{komatsu2016direct}, fluid-structure interactions~\cite{engels2015numerical} and thermal flows~\cite{cui2018coupled}. So far, IBM research for high-order methods has been relatively unexplored, with efforts focused on Poisson problems \cite{lew2008,lew2011} and cut-cell approaches \cite{fidkowski2007triangular,muller2017high}. In the context of volume penalization, we have recently extended this approach to high-order flux reconstruction schemes \cite{kou2021IBMFR1,kou2021IBMFR2}, and now to the high-order discontinuous Galerkin Spectral Element Method (DGSEM) in this work.

There have been several attempts to analyze the errors of the IBM approach. Bever and Leveque \cite{beyer1992analysis} analyzed the error of traditional IBM applied to one-dimensional problems, and highlighted the importance of choosing appropriate discrete delta functions to maintain optimal accuracy. Following a similar strategy, the immersed interface method \cite{li2015convergence}, which modifies the finite difference scheme with a jump condition for the immersed boundary, was derived \cite{leveque1994immersed}. Tornberg and Engquist \cite{tornberg2004numerical} performed error analyzes of traditional IBM with regularization and found first-order convergence for the standard central difference scheme with smoothing discrete delta functions. This error analysis was then extended to Stokes flows by Mori \cite{mori2008convergence}, Chen et al. \cite{chen2011note}, and Liu and Mori \cite{liu2014p} where error estimates for velocity and pressure were reported. Most analyses focus on the traditional IBM method, where the numerical property is based on the selection of the appropriate delta functions. Error analyzes for the direct forcing approach were also explored in \cite{guy2010accuracy,zhou2021analysis}, where the importance of maintaining smoothness in the solution and the choice of a suitable temporal and spatial resolution was highlighted. For these types of approach, the discretization error from space-time discretization and the modeling error from particular IBM treatment are coupled. In contrast, volume penalization has the advantage that the modeling error and the numerical error can be handled separately. The convergence of modeling errors was studied rigorously by Angot et al.~\cite{angot1999penalization} and Carbou and Fabrie~\cite{carbou2003boundary}. The modeling errors for the incompressible flow past a cylinder and a sphere were analyzed by Zhang and Zheng \cite{zhang2017high}. Therefore, the main concern is the discretization error, which depends on the numerical scheme used, and where a detailed error analysis is lacking, especially in the context of high-order methods.

In the present work, we perform error analyses of the IBM based on combination of volume penalization and nodal DGSEM, and propose new penalties to cancel spatial errors to improve the accuracy of the solution. In particular, we use modified equation analysis \cite{shyy1985errors,moura2015modified}, which has been extensively used to analyze the stability and accuracy of low-order numerical discretization, and to obtain high-order schemes \cite{warming1974modified,shubin1987modified}. The relationship between the errors introduced by the IBM based on volume penalization and high-order schemes remains unclear and motivates this work. First, using a modified equation analysis for volume penalization using DGSEM, we determine the shape and relationship of the dominant errors (i.e. dissipative/dispersive character). Second, we design the volume penalization scheme by including additional penalty terms that cancel the undesired numerical errors. In recent work, we have attempted to damp the numerical errors that arise from the volume penalization approach, using second-order derivatives \cite{kou2021VonNeumann} or combining it with a frequency damping technique \cite{kou2021SFD}. These studies have tried to minimize errors, without explicitly considering the causes of such errors, and therefore can be considered \textit{a posteriori} palliative treatment. In this work, we consider a different perspective and analyze the source of these errors. By doing so, we are able to cancel the errors at the source. This approach can be considered an \textit{a priori} error control. Note that we limit our analysis to linear advection-diffusion equations. The results from our analysis can thus be extended to linear systems (or linearized version of nonlinear equations). Examples include acoustics \cite{seo2006linearized} or stability analysis \cite{sipp_stability}. 

The article is organized as follows. In Section \ref{sec:motivation}, the volume penalization method and the DGSEM technique are introduced for the governing equation. Next, Section \ref{sec:EVP} introduces the principal errors in a volume penalization approach to investigate the discretization error using the modified equation analysis in Section \ref{sec:ME}. The numerical results are shown in Section \ref{sec:NumResults}, to validate the conclusions of the analysis, where one and two-dimensional advection-diffusion equations are investigated. Finally, conclusions are given in Section \ref{sec:Conclusion}. 

\section{Motivation}\label{sec:motivation}

\subsection{The governing equation}\label{subsec:govequation}
Let us introduce the problem by considering the following time-dependent 1D advection-diffusion equation for the transported solution $u = u(x,t)$,

\begin{subequations}\label{eqn:TransportProblem}
\begin{equation}
\begin{split}
\frac{\partial u}{\partial t} + c \frac{\partial u}{\partial x} - \nu \frac{\partial^2 u}{\partial x^2} = 0,\quad\text{for }x\in(0,L),\,t>0, 
\end{split}    
\end{equation}
where the flow parameters are constant: velocity field $c$ and kinetic viscosity $\nu \geq \nu_{\text{min}} > 0$. The PDE is completed with the set of initial and boundary conditions,

\begin{align}
& u(x,0) = u^{0}(x), & & 0\leq x\leq L, \\
& u(0,t) = u_{0}(t),\quad u(L,t) = u_{L}(t), & & t\geq 0.
\end{align} 
\end{subequations}
The transport problem (\ref{eqn:TransportProblem}) is discretized in space based on a high-order DG method; in time, a Runge-Kutta method; and some of the solution points would be penalized by additional source terms to impose the IBM conditions.

\subsection{The volume penalization approach}\label{subsec:VP}
Motivated by the characteristic-based VP \cite{brown2014CBVP} and the inclusion of local dissipation \cite{kou2021VonNeumann}, we consider the governing equation with penalization terms for the solution (classic volume penalization) and additional first-order and second-order penalization terms:

\begin{subequations}
\begin{equation}
\begin{split}\label{eqn:CDeqaution}
\frac{\partial u}{\partial t} + c \frac{\partial u}{\partial x} - \nu \frac{\partial^2 u}{\partial x^2} + \frac{\chi}{\eta_{1}}\left(u - u_{\text{s}}\right) + \frac{\partial }{\partial x}\left(\frac{\chi}{\eta_{2}}u \right) + \frac{\partial^2}{\partial x^2}\left(\frac{\chi}{\eta_{3}} u \right) = 0, 
\end{split}    
\end{equation}
The additional term in Equation (\ref{eqn:CDeqaution}) helps to impose the IBM on a given region of the domain $\Omega = [0,L]$. In this work, we consider a boundary condition of homogeneous Dirichlet type, namely $u_{\text{s}} = u_{\text{s}}(x,t) = 0$ (that is, a non-slip wall). The other two parameters are penalized terms determined in \textit{the modified equation analysis} section, where we focus only on the spatial errors of the discretization. The penalization parameters for variable, first and second order derivatives are $\eta_1$, $\eta_2$, and $\eta_3$ respectively; and to facilitate the analysis, a continuous mask function represented by a hyperbolic tangent function is defined as

\begin{align}
\chi = \chi(x,t) = \begin{dcases}
                   \left[\text{tanh}\left( d/\delta \right) + 1\right] /2, & \text{If } x\in\Omega_{\text{s}}\\
                   \left[\text{tanh}\left( -d/\delta \right) + 1\right] /2, & \text{Otherwise}
                   \end{dcases},
\end{align}
\end{subequations}
which distinguishes between the fluid, $\Omega_{\text{f}}$, and the solid, $\Omega_{\text{s}}$, regions such that $\Omega = \Omega_{\text{f}} \cup \Omega_{\text{s}}$. The distance of any solution point from the boundary interface is defined as $d = d(x,t)$. The width of the hyperbolic tangent function is defined as $\delta$, which should be infinitely small to reduce the modeling error and approximate the sharp mask function. This smooth mask ensures that spatial derivatives can be calculated. Note that for standard volume penalization, the mask function is sharp, which is 1 in the solid and 0 in the fluid. Later in Section \ref{subsec:1DTest}, we will show that both sharp and smooth masks result in similar behavior given a sufficiently small $\delta$.

\subsection{The VP-DG discrete equation}\label{subsec:DGSEM}
Equation (\ref{eqn:CDeqaution}) is discretized using the DG spectral element method. We group the terms in Equation (\ref{eqn:CDeqaution}) that leads to the penalized equation:
\begin{subequations}
\begin{align}\label{eqn:PDE}
\frac{\partial u}{\partial t} + \mathcal{L}u = 0,
\end{align}
where the second-order differential operator is represented by

\begin{align}\label{eqn:VPconvectiveflux}
\mathcal{L}u = \frac{\partial}{\partial x}\left(\widehat{c}u - \widehat{\nu}\frac{\partial u}{\partial x}\right) + \frac{\chi}{\eta_1}u,  
\end{align}
\end{subequations}
where $\widehat{c}$ and $\widehat{\nu}$ are the VP velocity field and the VP viscosity, respectively. Note that here we consider the element to be either a fully solid or a fully fluid one. This implies that the solid boundary aligns with the element interface, which is a natural choice for the DG method when using the local r-refinement, e.g. \cite{marcon2020rp}.

\begin{figure}[htbp!]
\centering
\includegraphics[width=1\textwidth]{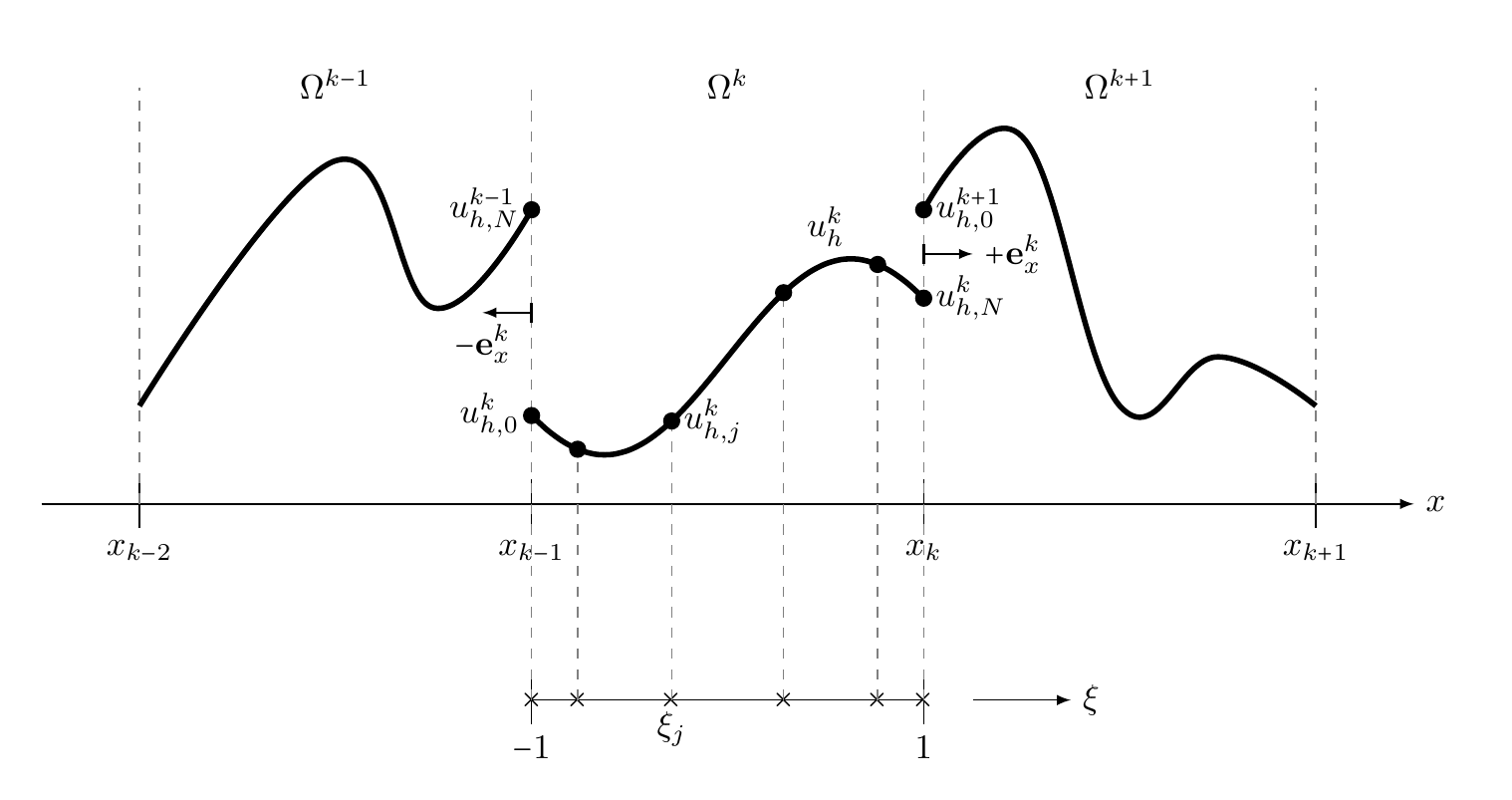}
\caption{Domain decomposition and reference interval in the DGSEM technique.}
\label{fig:DGSEM}
\end{figure}

The domain $\Omega$ is divided into multiple subdomains named elements $\Omega^k = [x_{k-1},x_{k}]$, $k = 1,2,\ldots,K$, as can be seen in Figure \ref{fig:DGSEM}, and mapped to the reference interval $\xi\in[-1,1]$. The global solution is assumed to be approximated by a piecewise polynomial defined as the direct sum $\oplus$ of the $K$ local polynomial solutions,

\begin{align}
u(x,t) \simeq u_h(x,t) = \bigoplus_{k=1}^Ku^k_h(x(\xi),t),
\end{align}
also for the test function and the VP flux function. On each element, we describe the local solution by the Legendre orthogonal interpolating polynomial, which is written in the Lagrange form, 

\begin{align}
u^k_h (\xi, t) = \sum_{j=0}^{N}u^k_{h,j}(t)\,l_j(\xi).
\end{align}
We select Gauss-Lobatto (GL) points, as they are becoming very popular in newly energy-stable and entropy-conserving schemes (\cite{manzanero2020entropystable,Wintermeyer2017entropystable}). Using GL points, the nodal (grid point) values become $u^k_{h,j}(t) = u^k_h(\xi_j,t)$, and $l_j(\xi)$ is the $N$-th order Lagrange interpolating polynomial, 

\begin{align}
l_j(\xi) \coloneqq \prod_{\substack{i=0 \\ i\neq j}}^{N}\frac{\xi - \xi_i}{\xi_j - \xi_i},
\end{align}
which satisfy $l_j(\xi_i) = \delta_{ij}$, being $\delta_{ij}$ the Kronecker delta. After obtaining weak forms of Equation (\ref{eqn:PDE}) and applying the Gaussian quadrature to the inner product in the reference interval (see Appendix \ref{appx:DGSEM} for details), the semi-discrete equation writes as follows: 

\begin{subequations}\label{eqn:DGSEM}
\begin{align}\label{eqn:semi_discrete}
\frac{\diff u^k_{h,j}}{\diff t} + \sum_{i=0}^{N}D^k_{ij}u^k_{h,i} = \numS^{k}_{j},
\end{align}
for $k = 1,2,\ldots,K$ and $j = 0,1,\ldots,N$ where

\begin{align}\label{eqn:D_matrix}
D^k_{ij} \coloneqq \frac{\chi^k_j}{\eta_1}\delta_{ij} - \frac{2}{\Delta x_{k}}\widehat{c}^{k}_{i}\frac{w_i}{w_j}l_j'(\xi_i) - \left(\frac{2}{\Delta x_{k}}\right)^{2}\widehat{\nu}^{k}_{i}\frac{w_i}{w_j}\sum_{r=0}^{N}l_j'(\xi_r)l_r'(\xi_i),
\end{align}
is the VP-DG derivative, and

\begin{align}\label{eqn:Num_Source}
\numS^{k}_{j} \coloneqq \frac{2}{\Delta x_k}\frac{\flux_{-1}^kl_j(-1) - \flux_1^kl_j(1)}{w_j} + \left(\frac{2}{\Delta x_k}\right)^2\frac{\widehat{\nu}^{k}_{0}\numU_{-1}^kl'_j(-1) - \widehat{\nu}^{k}_{N}\numU_1^kl'_j(1)}{w_j},
\end{align}
\end{subequations}
the numerical source. The weights of the GL quadrature are $\{w_i\}_{i=0}^{N}$ and $l'$ represents the derivative of the Lagrange polynomial. In the previous formula, we define the numerical flux as $\flux$ and its numerical counterpart as $\numU$. These fluxes are given by:

\begin{subequations}\label{eqn:weightsfandg}
\begin{align}
& \flux^{k}_{-1} = \mathfrak{f}^{k-1}_{2}u_{h,2}^{k-1} + \mathfrak{f}^{k}_{0}u_{h,0}^{k}, & & \numU^{k}_{-1} = \mathfrak{g}^{k-1}_{2}u_{h,2}^{k-1} + \mathfrak{g}^{k}_{0}u_{h,0}^{k}, \\
& \flux^{k}_{1} = \mathfrak{f}^{k}_{2}u_{h,2}^{k} + \mathfrak{f}^{k+1}_{0}u_{h,0}^{k+1}, & & \numU^{k}_{1} = \mathfrak{g}^{k}_{2}u_{h,2}^{k} + \mathfrak{g}^{k+1}_{0}u_{h,0}^{k+1}.
\end{align}
\end{subequations}
The weights $\mathfrak{f}$ and $\mathfrak{g}$ depend on $\widehat{c}$, $\widehat{\nu}$, and some numerical parameters that determine the advective/diffusive flux scheme used; see Appendix \ref{appx:DGSEM} for details. To calculate the viscous flux, we have considered the Bassi-Rebay 1 (BR1) scheme \cite{bassi1997DGViscous} and the Local discontinuous Galerkin (LDG) scheme \cite{Cockburn1998LDG}.

\section{Errors in volume penalization}\label{sec:EVP}
Rigorous proofs of the convergence of modeling errors have been provided in previous work \cite{angot1999penalization,carbou2003boundary}, showing that the numerical error introduced from the penalization term can be controlled \textit{a priori} \cite{brown2014CBVP}. Analysis of volume penalization suggests that the two contributions to total error are modeling and discretization errors \cite{engels2015numerical}:

\begin{align}
\left \|  u^{\text{exact}}  - u_{\eta} ^ {\text{num.}} \right \| \leq \left \|  u^{\text{exact}}  - u_{\eta} \right \| + \left \|  u_{\eta} - u_{\eta} ^ {\text{num.}} \right \|,
\end{align}
where $u^{\text{exact}}$ is the exact solution of the governing equation, $u_{\eta} $ and $u_{\eta} ^ {\text{num.}}$ are the exact and numerical solutions of the penalized equation, respectively, and $\left \|  \cdot \right \|$ is the $L_p$ norm used to quantify the error. The modeling error depends on the penalization parameter \cite{schneider2015immersed}:

\begin{align}
\left \|  u^{\text{exact}}  - u_{\eta} \right \| \propto \eta_1 ^ {\alpha}.
\end{align}

This explains that the convergence of the solution to the exact solution requires the error norm to approach zero for small penalization parameter limit, i.e., 

\begin{align}
\lim_{\eta_1 \rightarrow 0} \left \|  u^{\text{exact}}  - u_{\eta} \right \| = 0.   
\end{align}
According to Angot et al. \cite{angot1999penalization} and Carbou and Fabrie \cite{carbou2003boundary}, the volume
penalization gives $\alpha = 1/2$, indicating that the penalization error has a decay rate of $\mathcal{O}(\sqrt{\eta_1})$ for Dirichlet boundary conditions. For the Neumann boundary conditions, a decay rate of $\mathcal{O}(\eta_1)$ can be obtained \cite{kolomenskiy2015analysis}.

The second part of the overall error is the discretization error, which refers to the error between the exact solution and the numerical solution of the penalized equation. Details are given in the next section.

\section{The modified equation analysis}\label{sec:ME}
When we discretize the penalized equation (\ref{eqn:PDE}) numerically, we translate it into a semi-discrete system (\ref{eqn:semi_discrete}) for each element. This scheme is an approximation of our original equation. A different view is that the discrete system is the solution of modified differential equations but with some extra terms. This equation is named the modified equation (or reduced PDE):

\begin{align}
\frac{\partial u_h}{\partial t} + \mathcal{L}_hu + HOT = s(u_h^\star)
\end{align}
and is obtained by expanding the solutions in Equation (\ref{eqn:semi_discrete}) with a Taylor series around a point in the mesh $x(\xi_j)$. We omit the superscript $k$. $HOT$ is the high-order term due to the Taylor series. $\mathcal{L}_h$ is the operator $\mathcal{L}$ applied at $x(\xi_j)$ and $s(u_h^\star)$ is a function of $u_h^\star$ that is the solution transported from the element $\Omega^k$ at the interfaces; see Figure \ref{fig:DGSEM}. The purpose of the modified equation is to obtain the truncation error, $TE$, which is defined as the difference between the original equation and the modified equation. With a consistent discretization and a stable numerical scheme, the discretization error or the $TE$ term writes as follows:

\begin{align}
TE = C_0\Delta x^p + C_1\Delta x^{p + 1} + C_2\Delta x^{p + 2} + C_3\Delta x^{p + 3} + \ldots = \mathcal{O}(\Delta x^p),
\end{align}
where $\Delta x$ is a geometric discretization parameter representative of the grid spacing, $p$ the order of accuracy of the numerical scheme, and $C_0,\,C_1,\,C_2,\,\ldots$ some constants that do not depend on the discretization. However, the discretization error is not only determined by the numerical scheme, but is also limited by the regularity of the solution, as pointed out by Schneider et al. \cite{kadoch2012volume,schneider2015immersed}. For high-order methods, with good regularity at the interface, the high-order convergence property can be recovered, that is, $\mathbcal{O}(\Delta x ^ {N+1})$. This has been shown in the recent work of Kou et al. \cite{kou2021IBMFR2}. Here, we further analyze discretization errors to control and reduce these errors and improve accuracy. 

Consider that finite-volume/difference methods are local by nature. In contrast, DGSEM is local within the whole domain but global within the element. Due to the non-local character of DG within each element, the solution depends on every point at the GL mesh. To perform the analysis, we center the solution at the same point of the source for each component of the discrete equation. Let us simplify the analysis to three GL points ($N = 2$), which are located at $\xi_0 = -1$, $\xi_1 = 0$ and $\xi_2 = 1$ with weights $w_0 = 1/3$, $w_1 = 4/3$ and $w_2 = 1/3$ respectively. The Lagrange polynomials are 

\begin{align}
l_0 = \frac{1}{2}\xi(\xi - 1),\quad l_1 = -(\xi + 1)(\xi - 1),\quad l_2 = \frac{1}{2}\xi(\xi + 1),
\end{align}
and the VP-DG matrix,

\begin{align}
\begin{split}
\begin{pmatrix}
D^k_{00} & D^k_{10} & D^k_{20} \\
D^k_{01} & D^k_{11} & D^k_{21} \\
D^k_{02} & D^k_{12} & D^k_{22} 
\end{pmatrix} = \frac{1}{\eta_1} \begin{pmatrix}
                                 \chi^k_{0} & 0          & 0          \\
                                 0          & \chi^k_{1} & 0          \\
                                 0          & 0          & \chi^k_{2} 
                                 \end{pmatrix} &  - \frac{2}{\Delta x_k} \begin{pmatrix}
                                                                       -\dfrac{3}{2}\widehat{c}^k_{0} & -2\widehat{c}^k_{1} &  \dfrac{1}{2}\widehat{c}^k_{2} \\
                                                                        \dfrac{1}{2}\widehat{c}^k_{0} &  0              & -\dfrac{1}{2}\widehat{c}^k_{2} \\
                                                                       -\dfrac{1}{2}\widehat{c}^k_{0} &  2\widehat{c}^k_{1} & \dfrac{3}{2}\widehat{c}^{k}_{2} 
                                                                       \end{pmatrix} - \\
                                                & \left(\frac{2}{\Delta x_k}\right)^{2} \begin{pmatrix}
                                                                                         \widehat{\nu}_{0}^k             &  4\widehat{\nu}_{1}^k &  \widehat{\nu}_{2}^k \\
                                                                                        -\dfrac{1}{2}\widehat{\nu}_{0}^k & -2\widehat{\nu}_{1}^k & -\dfrac{1}{2}\widehat{\nu}_{2}^k \\
                                                                                         \widehat{\nu}_{0}^k             &  4\widehat{\nu}_{1}^k &  \widehat{\nu}_{2}^k \\
                                                                                        \end{pmatrix}.
\end{split}
\end{align}
comming from Equation (\ref{eqn:D_matrix}). Due to the non-local character of DG within each element, the solution depends on every point at the GL mesh. To perform the analysis, we center the solution at the same point of the source for each component of the discrete equation. For example, the discrete equation for the $j = 0$ component is

\begin{subequations}
\begin{align}\label{eqn:j=0_DiscreteEquation}
\frac{\diff u^k_{h,0}}{\diff t} + D^k_{00}u^k_{h,0} + D^k_{10}u^k_{h,1} + D^k_{20}u^k_{h,2} = \numS^k_0. 
\end{align}
and the numerical source,

\begin{align}
\numS^k_0 = \frac{2}{\Delta x_k}\left[3\flux^{k}_{-1} - \frac{3}{\Delta x_k}\left(3\widehat{\nu}^k_{0}\numU_{-1}^k + \widehat{\nu}^k_{2}\numU_{1}^k\right)\right]. 
\end{align}
\end{subequations}
Expanding $u^k_{h,1}$ and $u^k_{h,2}$ around $u^k_{h,0}$, we have

\begin{align}
& u^k_{h,1} = u^k_{h,0} + \Delta\xi\left.\frac{\partial u^k_h}{\partial\xi}\right|_{\xi_0} + \frac{\Delta\xi^2}{2!}\left.\frac{\partial^2 u^k_h}{\partial\xi^2}\right|_{\xi_0} + \frac{\Delta\xi^3}{3!}\left.\frac{\partial^3 u^k_h}{\partial\xi^3}\right|_{\xi_0} + \ldots \\
& u^k_{h,2} = u^k_{h,0} + 2\Delta\xi\left.\frac{\partial u^k_h}{\partial\xi}\right|_{\xi_0} + \frac{2^2\Delta\xi^2}{2!}\left.\frac{\partial^2 u^k_h}{\partial\xi^2}\right|_{\xi_0} + \frac{2^3\Delta\xi^3}{3!}\left.\frac{\partial^3 u^k_h}{\partial\xi^3}\right|_{\xi_0} + \ldots
\end{align}
being $\Delta\xi = \xi_1 - \xi_0 = \xi_2 - \xi_1 = 1$, we get

\begin{equation}
\begin{split}
\frac{\partial u^k_{h,0}}{\partial t} +& \left(D^k_{00} + D^k_{10} + D^k_{20}\right)u^k_{h,0} + \\
& \sum_{m = 1}^{\infty} \left(D^k_{10} + 2^m \left( D^k_{20} + \frac{6}{\Delta x_k^2}\widehat{\nu}_2^k\mathfrak{g}_2^k\right)\right)\frac{\Delta\xi^m}{m!}\left.\frac{\partial^m u^k_h}{\partial\xi^m}\right|_{\xi_0} = \numS^k_0,
\end{split}\label{eqn:preReduced_PDE}
\end{equation}
where the numerical source now reads:

\begin{align}
\begin{split}
\numS^k_0 = \frac{2}{\Delta x_k}\Bigg[\Bigg.&\left(3\mathfrak{f}_{2}^{k-1} - \frac{9}{\Delta x_k}\widehat{\nu}_{0}^{k}\mathfrak{g}_{2}^{k-1}\right)u_{h,2}^{k-1} - \frac{3}{\Delta x_k}\widehat{\nu}_{2}^{k}\mathfrak{g}_{0}^{k+1}u_{h,0}^{k+1} + \\
& \left(3\mathfrak{f}_{0}^{k} - \frac{3}{\Delta x_k}\left(3\widehat{\nu}_{0}^{k}\mathfrak{g}_{0}^{k} + \widehat{\nu}_{2}^{k}\mathfrak{g}_{2}^{k}\right)\right)u_{h,0}^{k}\Bigg.\Bigg].   
\end{split}
\end{align}
Neither $u_{h,2}^{k-1}$ nor $u_{h,0}^{k+1}$ can be expanded using Taylor series due to the discontinuous nature of DG. The terms in brackets of Equation (\ref{eqn:preReduced_PDE}) simplify to

\begin{align}
D^k_{00} + D^k_{10} + D^k_{20} = \frac{\chi^k_0}{\eta_1} + \frac{3\widehat{c}^k_{0} + 4\widehat{c}_{1}^k - \widehat{c}_{2}^k}{\Delta x_k} - 4\frac{\widehat{\nu}_{0}^k + 4\widehat{\nu}_{1}^k + \widehat{\nu}_{2}^k}{\Delta x_k^2} \eqqcolon \frac{\chi^k_0}{\eta_1} + \widetilde{r}^k_{0},
\end{align}
and 

\begin{align}
\begin{split}
D^k_{10} + 2^m\left(D^k_{20} + \frac{6}{\Delta x_k^2}\widehat{\nu}_2^k\mathfrak{g}_2^k\right) = &3\frac{2^{m+1}}{\Delta x_k^2}\widehat{\nu}_2^k\mathfrak{g}_2^k + 4\frac{\widehat{c}^k_{1} - 2^{m-2}\widehat{c}^k_{2}}{\Delta x_k} - \\ & 16\frac{\widehat{\nu}^k_{1} + 2^{m-2}\widehat{\nu}_{2}^k}{\Delta x_k^2} \eqqcolon \left(\frac{2}{\Delta x_k}\right)^m\Zhe_{0}^{(m)k}.
\end{split}
\end{align}
For the Taylor term of first order, $\widetilde{c}^k_{0} = \Zhe_{0}^{(1)k}$; and the second order term, $\widetilde{\nu}^k_{0} = -\Zhe_{0}^{(2)k} / 2$. Finally, we find the modified equation at the left-boundary element,

\begin{subequations}\label{eqn:Reduced_PDE_LBE}
\begin{equation}
\begin{split}
\frac{\partial u^k_{h,0}}{\partial t} + \Delta\xi\,\widetilde{c}^k_{0}\frac{2}{\Delta x_k}\left.\frac{\partial u^k_h}{\partial\xi}\right|_{\xi_0} -& \Delta\xi^2\,\widetilde{\nu}^k_{0}\left(\frac{2}{\Delta x_k}\right)^2\left.\frac{\partial^2 u^k_h}{\partial\xi^2}\right|_{\xi_0} + \\
                                    & \frac{\chi^k_0}{\eta_1}u^k_{h,0} + HOT_0^k = s^k_{\text{DG},0},\label{eqn:3pointGL_Reduced_PDE}
\end{split}
\end{equation}
where 

\begin{align}
s^k_{\text{DG},0} = \numS^k_0 - \widetilde{r}^k_{0}u^{k}_{h,0},\label{eqn:s_DG}
\end{align}
and

\begin{align}
HOT_0^k = \sum_{m = 3}^{\infty} \left(\frac{2}{\Delta x_k}\right)^m\Zhe_{0}^{(m)k}\frac{\Delta\xi^{m}}{m!}\left.\frac{\partial^m u^k_h}{\partial\xi^m}\right|_{\xi_0}.\label{eqn:HOT0_v1}
\end{align}
Additionally, the original PDE at the left-boundary element is 

\begin{align}
\frac{\partial u^k_{h,0}}{\partial t} + \widehat{c}_{0}^{k}\frac{2}{\Delta x_k}\left.\frac{\partial u^k_{h}}{\partial \xi}\right|_{\xi_{0}} - \widehat{\nu}_{0}^k\left(\frac{2}{\Delta x_k}\right)^2\left.\frac{\partial^2 u_h^k}{\partial \xi^2}\right|_{\xi_{0}} + \frac{\chi_{0}^k}{\eta_{1}}u^k_{h,0} = 0.
\end{align}
and, therefore, the truncation error at the left-boundary element becomes: 

\begin{equation}
\begin{split}
TE_0^k = s^k_{\text{DG},0} +& \left(\widehat{c}_{0}^{k} - \Delta\xi\,\widetilde{c}^k_{0} \right)\frac{2}{\Delta x_k}\left.\frac{\partial u^k_{h}}{\partial \xi}\right|_{\xi_{0}} - \\
                            & \left(\widehat{\nu}_{0}^k - \Delta\xi^2\,\widetilde{\nu}^k_{0}\right)\left(\frac{2}{\Delta x_k}\right)^2\left.\frac{\partial^2 u_h^k}{\partial \xi^2}\right|_{\xi_{0}} - HOT_{0}^{k}.\label{eqn:3pointGL_TE}
\end{split}
\end{equation}
\end{subequations}

We can proceed in a similar manner to obtain the modified equations and truncation errors for the inner point, $j = 1$, and the right-boundary point, $j = 2$. For $j = 1$, $u^k_{h,0}$ and $u^k_{h,2}$ are centered on $u^k_{h,1}$; for $j = 2$, $u^k_{h,0}$ and $u^k_{h,1}$ are centered on $u^k_{h,2}$. Their formulae can be written using equations (\ref{eqn:Reduced_PDE_LBE}), but with differences in the reactive parameter, $\widetilde{r}^k_j$, the coefficient $\Zhe$ and the numerical source, $\numS_j^k$, for $j = 0,1,2$, see Table \ref{tab:j_Reduced_PDE} and \ref{tab:NumericalSource}. The source of the DG, $s^k_{\text{DG},j}$, arises from the discontinuous nature of the DG approach (discontinuous boundary values) and the selected diffusive scheme. 

{
\renewcommand{\arraystretch}{2}
\begin{table}[!htpb]
\centering
% \begin{adjustbox}{angle=0}
 \resizebox{\textwidth}{!}{
\begin{tabular}{cccc}
\hline
$j$ & $\xi_j$ & $\widetilde{r}_{j}^{k}$ & $\Zhe_{j}^{(m)k}$        \\ \hline
$0$ & $-1$    & $\dfrac{3\widehat{c}^k_{0} + 4\widehat{c}_{1}^k - \widehat{c}_{2}^k}{\Delta x_k} - 4\dfrac{\widehat{\nu}_{0}^k + 4\widehat{\nu}_{1}^k + \widehat{\nu}_{2}^k}{\Delta x_k^2}$                     & $\dfrac{2^{2-m}\widehat{c}^k_{1} - \widehat{c}^k_{2}}{\Delta x_k^{1-m}} - 4\dfrac{2^{2-m}\widehat{\nu}^k_{1} + \widehat{\nu}_{2}^k\left(1-\sfrac{3}{2}\mathfrak{g}_{2}^{k}\right)}{\Delta x_k^{2-m}}$        \\ 
$1$ & $0$     & $-\dfrac{\widehat{c}^k_{0} - \widehat{c}_{2}^k}{\Delta x_k} + 2\dfrac{\widehat{\nu}_{0}^k + 4\widehat{\nu}_{1}^k + \widehat{\nu}_{2}^k}{\Delta x_k^2}$                     & $- 2^{-m}\dfrac{(-1)^{m}\widehat{c}^k_{0} - \widehat{c}^k_{2}}{\Delta x_k^{1-m}} + 2^{1-m}\dfrac{(-1)^{m}\widehat{\nu}^k_{0}\left(1+3\mathfrak{g}_{0}^{k}\right) + \widehat{\nu}_{2}^k\left(1+3\mathfrak{g}_{2}^{k}\right)}{\Delta x_k^{2-m}}$ \\ 
$2$ & $1$     & $\dfrac{\widehat{c}^k_{0} - 4\widehat{c}_{1}^k - 3\widehat{c}_{2}^k}{\Delta x_k} - 4\dfrac{\widehat{\nu}_{0}^k + 4\widehat{\nu}_{1}^k + \widehat{\nu}_{2}^k}{\Delta x_k^2}$                     & $(-1)^{m}\left(\dfrac{\widehat{c}^k_{0} - 2^{2-m}\widehat{c}^k_{1}}{\Delta x_k^{1-m}} - 4\dfrac{\widehat{\nu}^k_{0}\left(1-\sfrac{3}{2}\mathfrak{g}_{2}^{k}\right) + 2^{2-m}\widehat{\nu}_{1}^k}{\Delta x_k^{2-m}}\right)$   \\ \hline
\end{tabular}
 }
% \end{adjustbox}
\caption{The reaction parameter and the coefficient $\Zhe$ in the modified equations for a three-point GL grid.}
\label{tab:j_Reduced_PDE}
\end{table}
}

{
\renewcommand{\arraystretch}{2}
\begin{table}[!htpb]
\centering
% \begin{adjustbox}{angle=0}
 \resizebox{\textwidth}{!}{
\begin{tabular}{cccc}
\hline
$j$ & $\numS_{j}^{k}$ \\ \hline
$0$ & $\dfrac{2}{\Delta x_k}\left[\left(3\mathfrak{f}_{2}^{k-1} - \dfrac{9}{\Delta x_k}\widehat{\nu}_{0}^{k}\mathfrak{g}_{2}^{k-1}\right)u_{h,2}^{k-1} - \dfrac{3}{\Delta x_k}\widehat{\nu}_{2}^{k}\mathfrak{g}_{0}^{k+1}u_{h,0}^{k+1} + \left(3\mathfrak{f}_{0}^{k} - \dfrac{3}{\Delta x_k}\left(3\widehat{\nu}_{0}^{k}\mathfrak{g}_{0}^{k} + \widehat{\nu}_{2}^{k}\mathfrak{g}_{2}^{k}\right)\right)u_{h,0}^{k}\right]$              \\
$1$ & $\dfrac{6}{\Delta x_k^2}\left[\widehat{\nu}_0^k\mathfrak{g}_2^{k-1}u_{h,2}^{k-1} + \widehat{\nu}_2^k\mathfrak{g}_0^{k+1}u_{h,0}^{k+1} + \left(\widehat{\nu}_{0}^k\mathfrak{g}_0^k + \widehat{\nu}_{2}^k\mathfrak{g}_2^k\right)u_{h,1}^k\right]$              \\
$2$ & $-\dfrac{2}{\Delta x_k}\left[\dfrac{3}{\Delta x_k}\widehat{\nu}_{0}^{k}\mathfrak{g}_{2}^{k-1}u_{h,2}^{k-1} + \left(3\mathfrak{f}_{0}^{k+1} + \dfrac{9}{\Delta x_k}\widehat{\nu}_{2}^{k}\mathfrak{g}_{0}^{k+1}\right)u_{h,0}^{k+1} + \left(3\mathfrak{f}_{2}^{k} + \dfrac{3}{\Delta x_k}\left(3\widehat{\nu}_{2}^{k}\mathfrak{g}_{2}^{k} + \widehat{\nu}_{0}^{k}\mathfrak{g}_{0}^{k}\right)\right)u_{h,2}^{k}\right]$              \\ \hline
\end{tabular}
 }
% \end{adjustbox}
\caption{Numerical source for the DG source, $s^k_{\text{DG},j} = \numS^k_j - \widetilde{r}^k_{j}u^{k}_{h,j}$.}
\label{tab:NumericalSource}
\end{table}
}

Now suppose that an element, $\Omega^k$, belongs to a solid region, $\Omega_{\text{s}}$, then $\chi_j^k = 1$ for $j=0,1,2$, and the truncation error still remains inside the solid. If we want to eliminate all the error terms in $TE_j^k$ for $j = 0,1,2$, we need to solve the system:

\begin{align}\label{eqn:KillingSystem}
\left\{\begin{array}{l}
       s^k_{\text{DG},j} = 0, \\
       \widehat{c}_{j}^{k} - \Delta\xi\,\widetilde{c}^k_{j} = 0, \\
       \widehat{\nu}_{j}^{k} - \Delta\xi\,\widetilde{\nu}^k_{j} = 0, \\
       \Zhe_{j}^{(m)k} = 0,
       \end{array}\right.    
\end{align}
for $j=0,1,2$ and $m\geq 3$. However, the problem is given by Eq. (\ref{eqn:KillingSystem}) has an infinite number of equations and finite unknowns. In total, there are 10 unknowns, which are the 4 weights $\mathfrak{f}$ and $\mathfrak{g}$, $\widehat{c}_0^k = \widehat{c}_1^k = \widehat{c}_2^k = c + 1/\eta_2 \eqqcolon \widehat{c}$ and $\widehat{\nu}_0^k = \widehat{\nu}_1^k = \widehat{\nu}_2^k = \nu - 1/\eta_3 \eqqcolon \widehat{\nu}$. The solution to the system is as follows.

\begin{subequations}
\begin{align}
\left\{\begin{array}{l}
\eta_2 = -1/c,\quad \eta_3 = 1/\nu \\
\mathfrak{f}_{2}^{k-1} = \mathfrak{f}_{0}^{k} = \mathfrak{f}_{2}^{k} = \mathfrak{f}_{0}^{k+1} = 0 \\
\text{for all } \mathfrak{g}_{2}^{k-1},\mathfrak{g}_{0}^{k},\mathfrak{g}_{2}^{k}\text{ and }\mathfrak{g}_{0}^{k+1}
\end{array}\right.
\end{align}
\end{subequations}
This solution will be referred to as the trivial solution of the problem. At this point, one may wonder if there is any other set, a nontrivial family, that cleans up almost all the errors within the solid region. To investigate this, a determined system should be formed. Ideal errors remaining within the solid region should be:

\begin{align}
TE_j^k \sim HOT_j^k \sim \mathcal{O}(\Delta x_k^m),\quad j = 0,1,2,    
\end{align}
for $m \geq 3$ as a representation of high order. However, the investigation of non-trivial solutions did not meet the previous requirement; see more details in the Appendix \ref{appx:Nontrivial}. Table \ref{tab:nontrivial} summarizes both the trivial and non-trivial solutions that have been found.

The trivial solution is the condition for DGSEM to compensate (or kill) spatial truncation errors within the body region, but additional insight can be obtained. If we substitute the values for $\eta_2$ and $\eta_3$ in our penalized equation and isolate $\chi$, we get the following:

\begin{align}
\frac{\partial u}{\partial t} + \underbrace{\frac{\partial}{\partial x}\left[(1-\chi)\left(cu - \nu\frac{\partial u}{\partial x}\right)\right]}_{\text{Physical term}} + \frac{\chi}{\eta_{1}}u = 0.
\end{align}
If we are in the solid region, $\chi = 1$, then the physical contribution of the PDE is removed, so this region is modeled with only the reaction penalization term, and therefore only time integration methods will lead to errors within the solid. This result agrees with the use of a typical characteristic-based volume penalization approach \cite{brown2014CBVP,abalakin2016immersed}, where the RHS vanishes to smooth out the errors in the solid region but without a theoretical explanation, which is provided here. Cancellation of particular terms can reduce the error inside the solid, thus improving the accuracy in the fluid region.

{
\renewcommand{\arraystretch}{2}
\begin{table}[!htpb]
\centering
 \resizebox{\textwidth}{!}{
\begin{tabular}{cccccc}
\hline
$\widehat{c} = c + \dfrac{1}{\eta_2}$ & $\widehat{\nu} = \nu - \dfrac{1}{\eta_3}$ & $\mathfrak{f}$s & $\mathfrak{g}$s & $\equiv$ CG * & $TE_j^k$ \\ \hline
$\eta_2$ free                         & $\eta_3$ free                             & free            & $\!\begin{aligned}[t]
& \mathfrak{g}^k_2 = \mathfrak{g}^k_0 = 2 \\ & \mathfrak{g}^{k-1}_2 = \mathfrak{g}^{k+1}_0 = 0
\end{aligned}$                  &            & $\mathcal{O}\left(\Delta x_k^0\right)$ \\  
$\eta_2$ free                         & $\eta_3$ free                             & $\!\begin{aligned}[t]
& \mathfrak{f}^{k-1}_2 = -\mathfrak{f}^k_2 = \widehat{c} + \frac{4}{\Delta x_k}\widehat{\nu} \\ & \mathfrak{f}^{k}_0 = \mathfrak{f}^{k+1}_0 = 0
\end{aligned}$ & $\!\begin{aligned}[t]
& \mathfrak{g}^k_2 = \mathfrak{g}^k_0 = 2 \\ & \mathfrak{g}^{k-1}_2 = \mathfrak{g}^{k+1}_0 = 0
\end{aligned}$                  & \Checkmark & $\mathcal{O}\left(\Delta x_k^0\right)$ \\ 
\multicolumn{2}{c}{$\widehat{c} + \dfrac{4}{\Delta x_k}\widehat{\nu} = 0$}                                      & 0 & $\!\begin{aligned}[t]
& \mathfrak{g}^k_2 = \mathfrak{g}^k_0 = 2 \\ & \mathfrak{g}^{k-1}_2 = \mathfrak{g}^{k+1}_0 = 0
\end{aligned}$                  &            & $\mathcal{O}\left(\Delta x_k^0\right)$ \\ 
$\eta_2$ free                         & 0                                         & $\!\begin{aligned}[t]
& \mathfrak{f}_{0}^{k} = -\mathfrak{f}_{2}^{k} = \widehat{c} \\ & \mathfrak{f}^{k}_0 = \mathfrak{f}^{k+1}_0 = 0
\end{aligned}$ & free                 & \Checkmark & $\mathcal{O}\left(\Delta x_k^2\right)$ \\ 
$\eta_2$ free                         & 0                                         & free & free                 &  & $\!\begin{aligned}[t]
& \mathcal{O}\left(\Delta x_k^0\right) \text{ Boundary} \\ & \mathcal{O}\left(\Delta x_k^2\right) \text{ Inner}
\end{aligned}$ \\ 
0 ($\eta_2 = -1/c$)                            & 0 ($\eta_3 = 1/ \nu$)                                                & 0 & free                 &  & $\mathcal{O}\left(\Delta x_k^\infty\right)$ \\ \hline
\end{tabular}
 }
\caption{Summary of family of solutions for VP-IBM DGSEM, the trivial solution is the last row. $*$ means equivalent to a continuous Galerkin (CG) method.}
\label{tab:nontrivial}
\end{table}
}

%\subsection{$N$-th order polynomials}\label{subsec:Nth-order}

%It can be seen that  and we set the condition (\ref{eqn:cond_VP}), the discrete equation becomes 
%\begin{align}
%\frac{\diff u^k_{h,j}}{\diff t} + \frac{1}{\eta_1}u^k_{h,j} = 0,    
%\end{align}
%for $j = 0,1,\ldots,N$. Since the PDE in the solid is $\diff u/\diff t + u/\eta_{1} = 0$, then $TE_{j}^{k} = 0$ $\forall j = 0,1,2,\ldots,N$ whatever $N$-th order polynomial.

%\textcolor{red}{WHAT IS THIS METHOD? THE NON-TRIVIAL CASE This can justify the use of a typical characteristic-based volume penalization \cite{brown2014CBVP,abalakin2016immersed} approach, where the RHS term vanishes inside the solid while only the penalization term remains. The cancellation of particular terms can reduce the error inside the solid, thus providing improved accuracy in the fluid region. }

Finally, let us note that the modified equation analysis, based on Taylor expansions, is more general than the one first anticipated for the considered advection-diffusion problem. Indeed, the Taylor analysis can be reframed into a more general Fourier analysis framework, as detailed in Appendix \ref{appx:FourierAndTaylor}, where we show that Fourier analysis links with the Taylor series. 

\section{Numerical results}\label{sec:NumResults}
In this section we introduce two numerical experiments to evaluate and validate the trivial solution derived from the modified equation analysis. The first group of cases is the one-dimensional advection-diffusion equation, where the influence of penalization parameters is studied in detail. The optimal parameters obtained from the numerical experiments and the analysis of the modified equations are then applied to the two-dimensional advection-diffusion equation.

\subsection{One-dimensional advection-diffusion equation}\label{subsec:1DTest}

\begin{figure*}[htbp]
		\centering
		\includegraphics[width=0.85\textwidth]{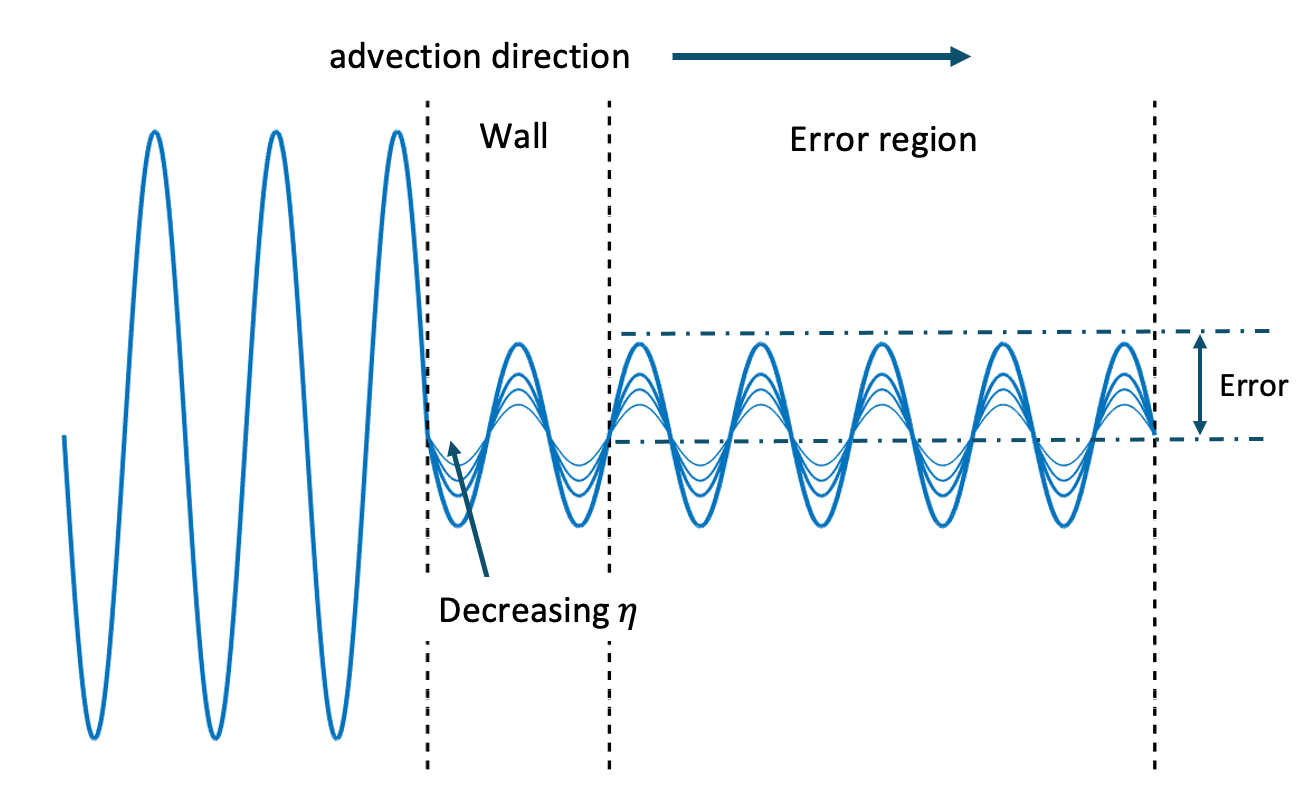}
		\caption{Schematic illustration of the advection problem with IBM.}
		\label{fig:cartoon1}
\end{figure*}

We start from the one-dimensional advection equation, where a non-slip wall is placed in the middle of the computational domain. This problem has been formulated in previous works \cite{kou2021VonNeumann,kou2021SFD}, which is illustrated in Figure \ref{fig:cartoon1}. The advection speed is set to $c = 1$ and the computational domain is defined in $x \in [-1,1]$, discretized by $K$ equispaced elements with mesh size $\Delta x$. The solution points are selected according to the Gauss-Lobatto quadrature rule, which is consistent with the previous analysis. Periodic boundary conditions are imposed on both sides of the computational domain. An upwind flux for the advection term is selected. The solid region is defined as a non-slip wall, i.e., $u_{\text{s}} = 0$. It lies in the middle of the computational domain and starts from $x=0$, whose width is defined as $\Delta_{\text{s}}$, leading to the solid region $0 \leq x \leq \Delta_{\text{s}}$.

For consistency with the analysis of the modified equations, we consider $\Delta_{\text{s}} = \Delta x$, which means that the solid boundaries lie exactly at the interface between the elements (if we have an even number of elements). This allows us to define the solid ratio $r = 1/K$ as the ratio between the solid region and the computational domain. As shown in Figure \ref{fig:cartoon1}, the initial wavelike solution passes through the non-slip wall in the middle, and the damped solution moves to the right as time evolves. Since periodic boundary conditions are considered, the solution will eventually become $0$. If the non-slip wall boundary condition is exactly imposed, the solutions coming out of the wall will be zero. However, in practice, this solution depends on the damping provided by the volume penalization approach, where a smaller $\eta$ makes the solution closer to zero. Therefore, the accuracy of the IBM imposition can be evaluated by comparing the exact solution (zero) and the damped solution in both the flow and solid regions (e.g., within a short advection time $0 < x< t$ ). 

We first perform the numerical experiment of a linear advection equation with a wavelike initial condition. The initial condition is defined as a sinusoidal wave with wavenumber $\omega$, which is nondimensionalized by the mesh size $\Delta x$ and the polynomial order $N$, defined as $\omega \Delta x/(N+1)$. Furthermore, due to the existence of a solid wall, the actual fluid domain is shorter than the entire computational domain; therefore, the effective wavenumber in the fluid region is greater than $\omega$. This effective wavenumber is rescaled by the solid ratio $r$ to $ \overline{\omega} = \omega / (1-r)$ \cite{kou2021VonNeumann}. We consider a spatial discretization with $K = 40$ elements in the computational domain ($\Delta x=0.05$). Based on this mesh, we set $\Delta_{\text{s}} = \Delta x$ with $r = 1/40$ and choose $N = 3$ as a representative order for high-order methods. The initial condition with wavenumber $\overline{\omega} \Delta x /(N+1) = 0.3223$ is considered, which lies in the resolved wavenumber region of the scheme. The time integration is based on the third-order Runge-Kutta scheme. To reduce the temporal error, a sufficiently small time step is set to $\Delta t = 10^{-5}$. The final time is set to $1.1$ to obtain a sufficiently penalized solution in the right region of the computational domain.

\begin{figure*}[htbp]
    \begin{subfigure}{0.48\textwidth}
 		\centering
		\includegraphics[width=1\textwidth]{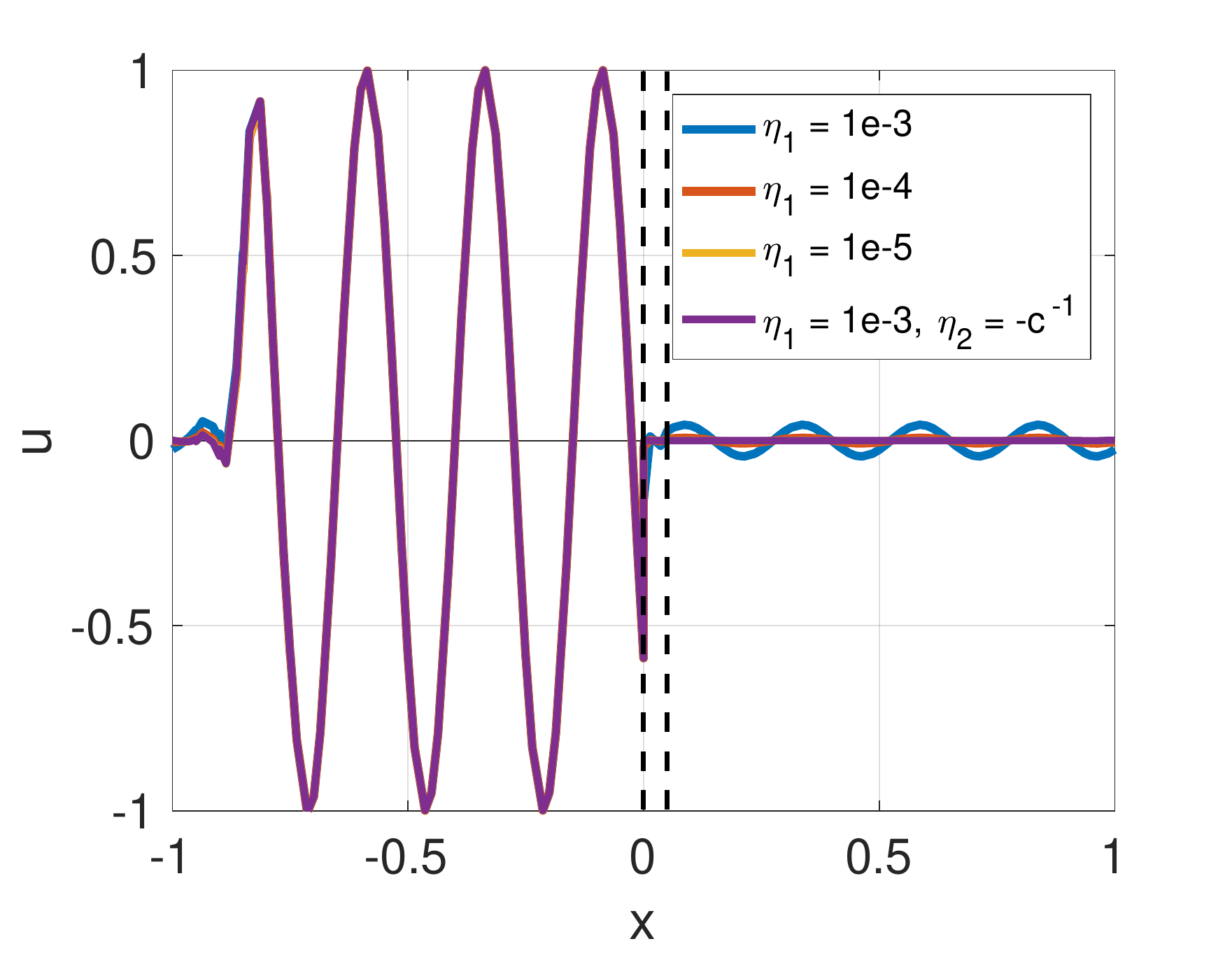}
		\caption{}
	\end{subfigure}
	\begin{subfigure}{0.48\textwidth}
 		\centering
		\includegraphics[width=1\textwidth]{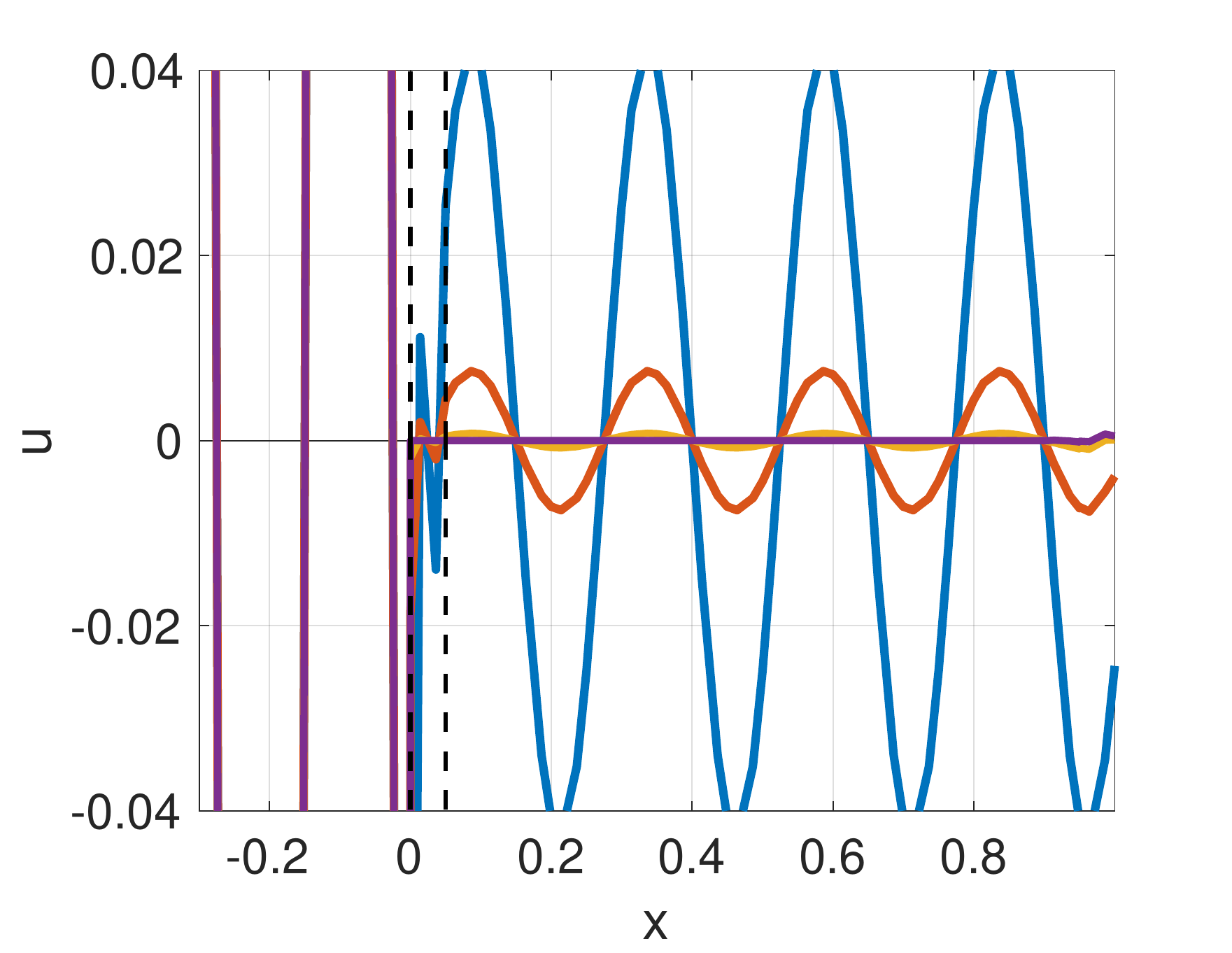}
		\caption{}
	\end{subfigure}
	\centering
	\caption{Simulation under different penalization parameters ($r = 1/40$, $N = 3$, initial wavenumber $\overline{\omega} \Delta x/(N+1) = 0.3223$, $K = 40$): a) Global view; b) Enlarged view.}
	\label{fig:advection-sim}
\end{figure*}

Different combinations of parameters (with and without the first-order term) are considered. To evaluate accuracy, the error (in the flow) is defined as the error in $x \in [\Delta_{\text{s}}, 1]$ and the penalized value $u_{\text{s}} = 0$. Defining the number of solution points inside the flow domain of interest as $N_{\text{p}} = (N+1)K$, we have the $L_2$-norm of the error as

\begin{align}
\text{error} = \sqrt{\frac{1}{N_{\text{p}}}\sum_{i=1}^{N_{\text{p}}} [u(x_i) - u^{\text{exact}}(x_i)]^2} \ , \ x_i \in [\Delta_{\text{s}}, 1], \ u^{\text{exact}} = 0,
\end{align}
and the $L_2$-norm of the error in the solid is defined as

\begin{align}
\text{error}_{\text{solid}} = \sqrt{\frac{1}{N_{\text{p}}}\sum_{i=1}^{N_{\text{p}}} [u(x_i) - u^{\text{exact}}(x_i)]^2} \ , \ x_i \in [0, \Delta_{\text{s}}], \ u^{\text{exact}} = 0,
\end{align}

\begin{figure}[htbp!]
    \includegraphics[width=1\textwidth]{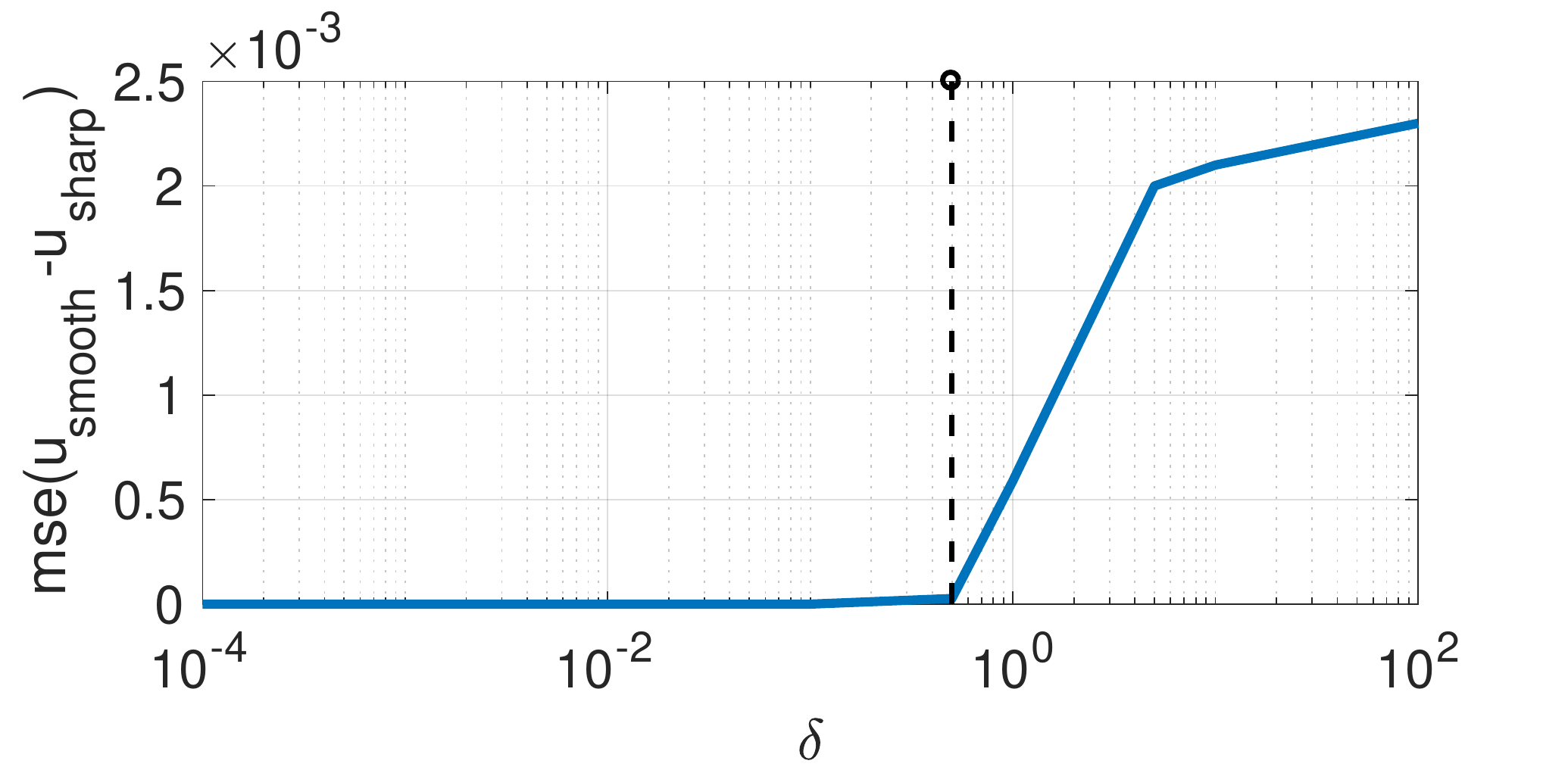}
	\centering
	\caption{Mean squared error between sharp and smooth mask function with increasing $\delta$.}
	\label{fig:mask}
\end{figure}

First, a numerical study is performed to justify the equivalence of using sharp or smooth mask functions (given the small width $\delta$ for the smooth mask function). Due to the fact that GL points overlap at the element interface, we can manually modify the smooth mask for the interface point corresponding to the solid element being $(\text{tanh}(1/\delta)+1)/2$ and the one in the fluid being $(\text{tanh}(-1/\delta)+1)/2$. We run the simulation at different $\delta$, $\eta_1 = 10^{-3}$, and final time 1.1, and compare the mean squared error between the results from the sharp and the smooth mask function. This error is compared in Figure \ref{fig:mask}, where the results based on the sharp or smooth mask are almost identical at small $\delta$, and the difference becomes dominant when $\delta$ is sufficiently large $\delta > 0.5$. Similar results are maintained when $\delta < 0.5$, which is sufficient to guarantee the equivalence of the analysis for smooth and sharp masks. As $\delta$ further increases, the mask becomes too smooth, resulting in additional modeling errors due to the wrong representation of the interface.

A comparison of the solution at the final time is shown in Figure \ref{fig:advection-sim}. Four cases are tested, where the first three cases contain only the volume penalization term for the solution, while the first-order penalization term with $\eta_2 = -1 / c$ is added to the last case. The figure shows that as the penalization parameter $\eta_1$ decreases, the solution approaches zero, indicating that the boundary condition is imposed more accurately. Note that in the last case, a large penalization parameter (i.e. weaker penalization) $\eta_1 = 10^{-3}$ is used. In addition, when the first-order term is added, improved accuracy is seen as the solution is closer to zero. The errors in the fluid region of the four cases are $3.071 \cdot 10^{-2}$, $5.385 \cdot 10^{-3}$, $5.698 \cdot 10^{-4}$, and $1.022 \cdot 10^{-4}$, respectively. This indicates that by introducing the first-order term with a proper selection of the penalization parameter, it is possible to improve the accuracy.

\begin{figure}[htbp!]
    \begin{subfigure}{.48\textwidth}
		\centering
		\includegraphics[width=1\textwidth]{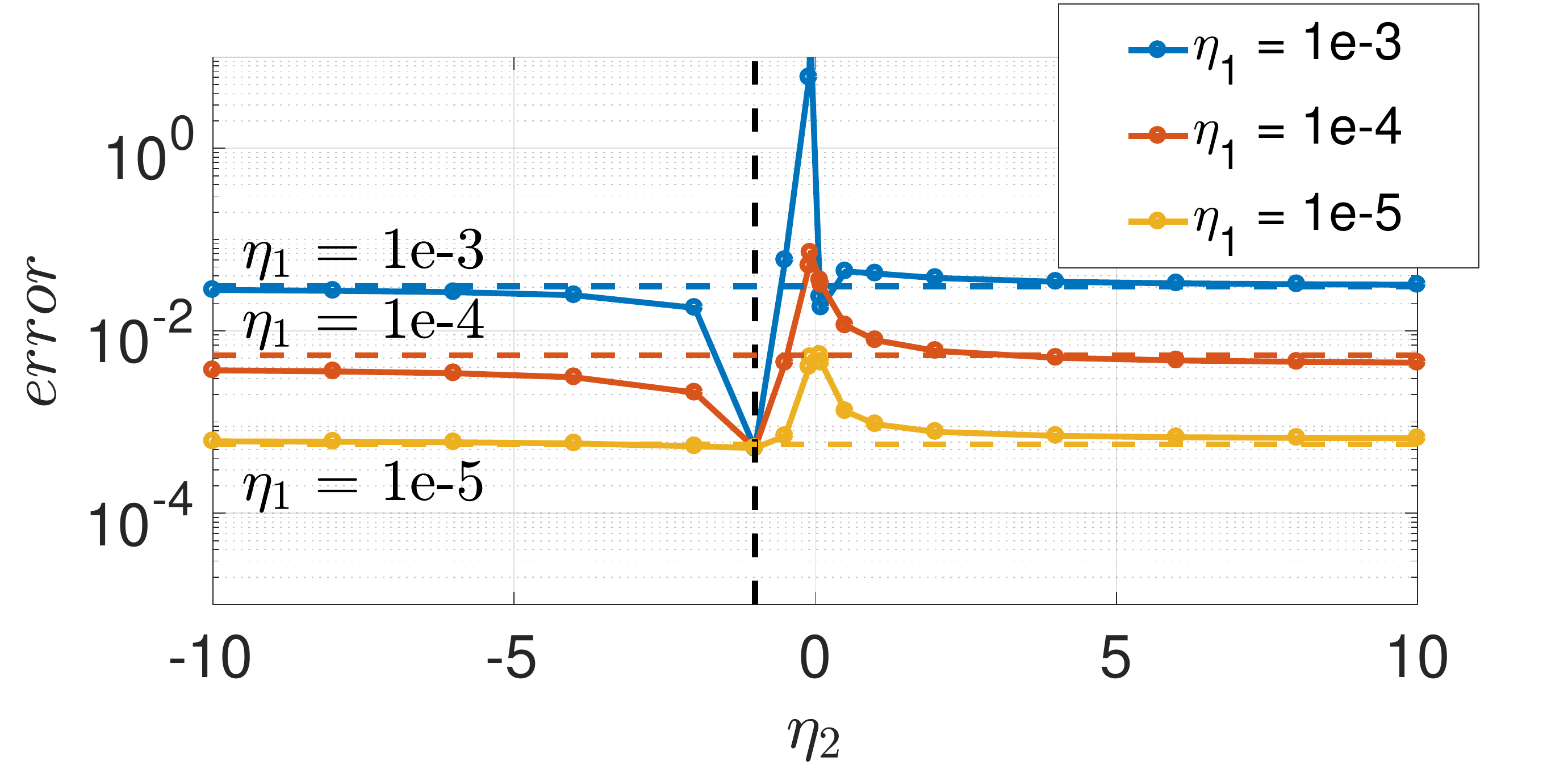}
		\caption{}
	\end{subfigure}
	\begin{subfigure}{.48\textwidth}
		\centering
		\includegraphics[width=1\textwidth]{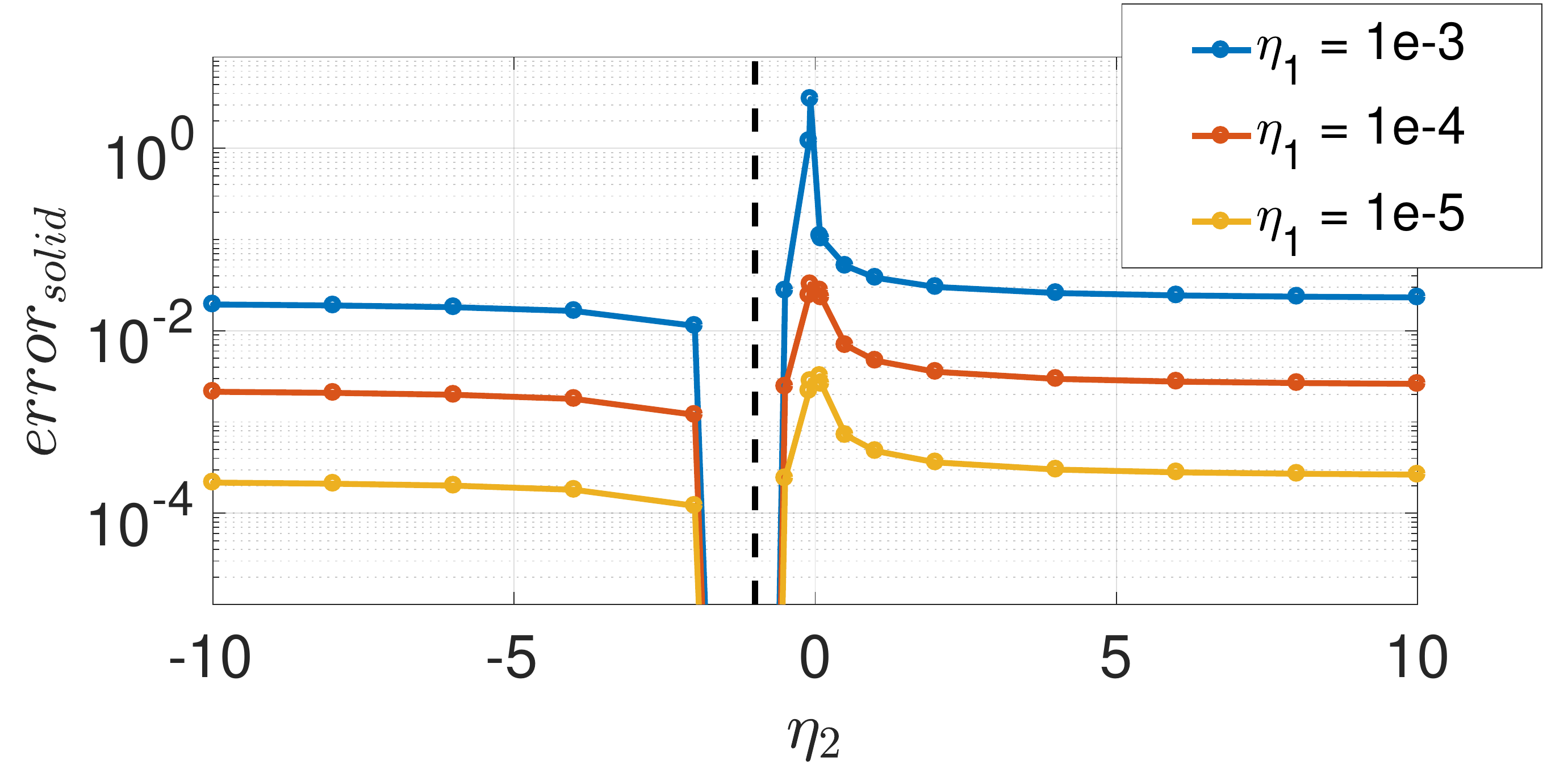}
		\caption{}
	\end{subfigure}
    \begin{subfigure}{.48\textwidth}
		\centering
		\includegraphics[width=1\textwidth]{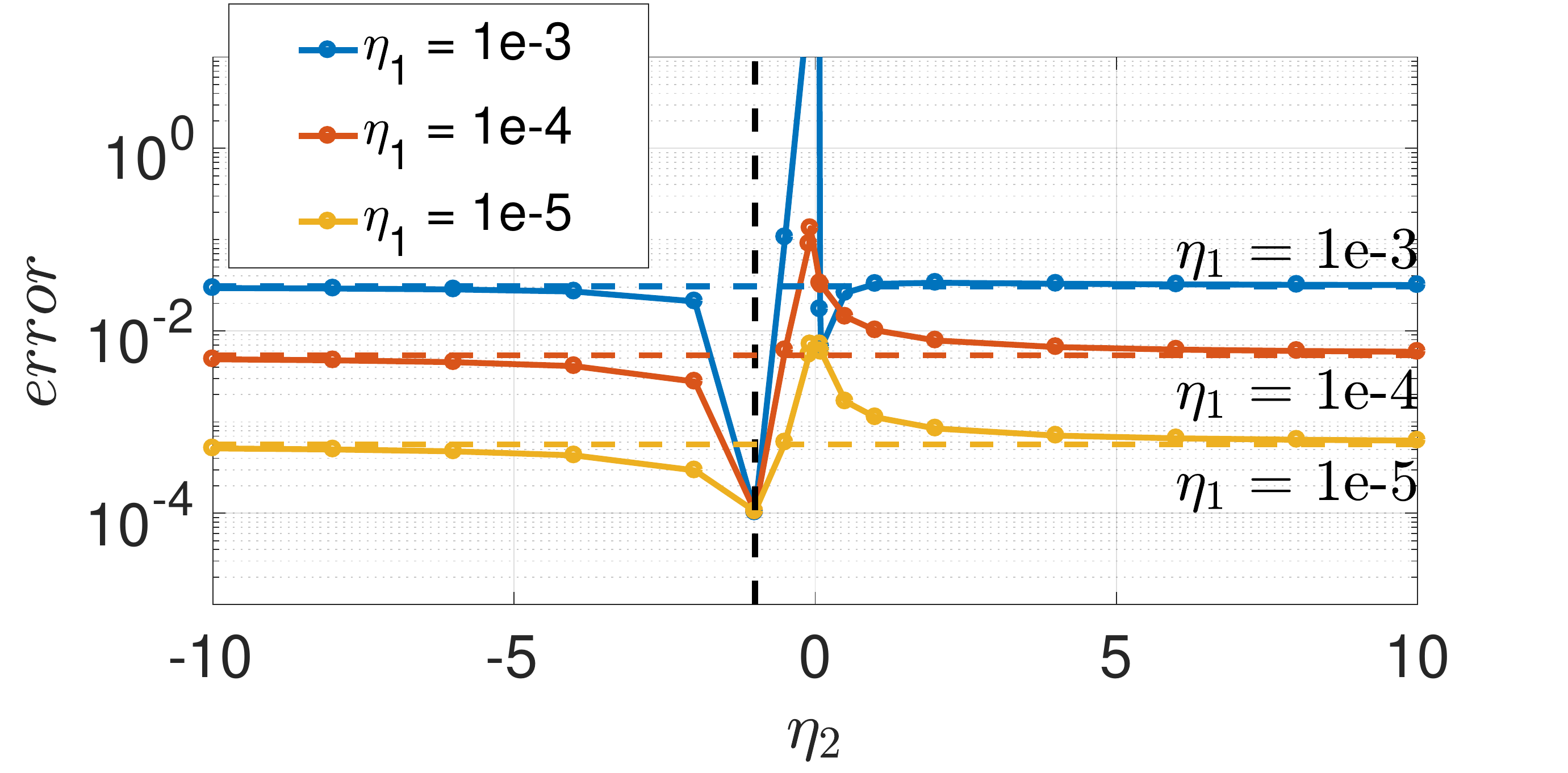}
		\caption{}
	\end{subfigure}
	\begin{subfigure}{.48\textwidth}
		\centering
		\includegraphics[width=1\textwidth]{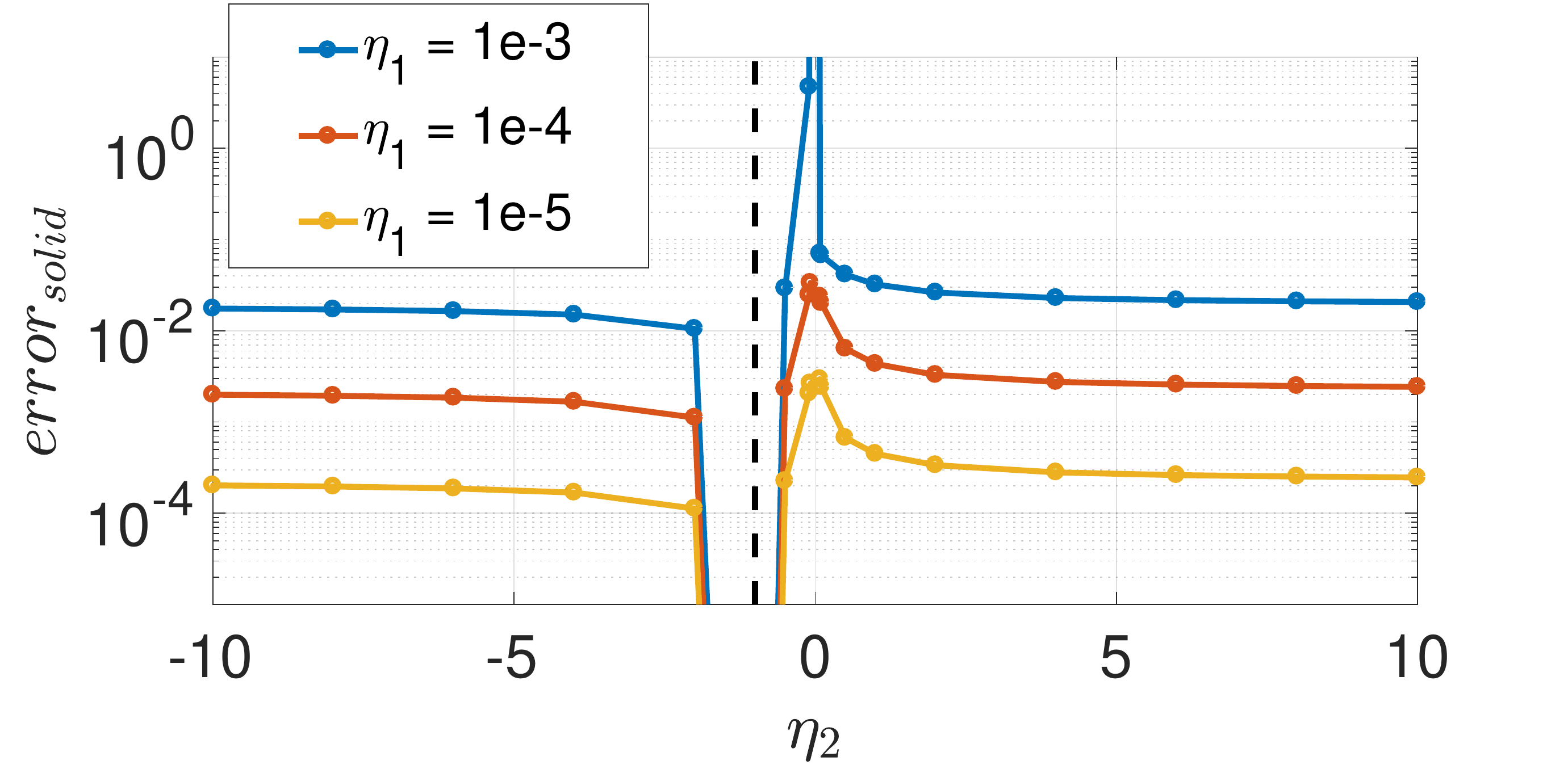}
		\caption{}
	\end{subfigure}
	\begin{subfigure}{.48\textwidth}
		\centering
		\includegraphics[width=1\textwidth]{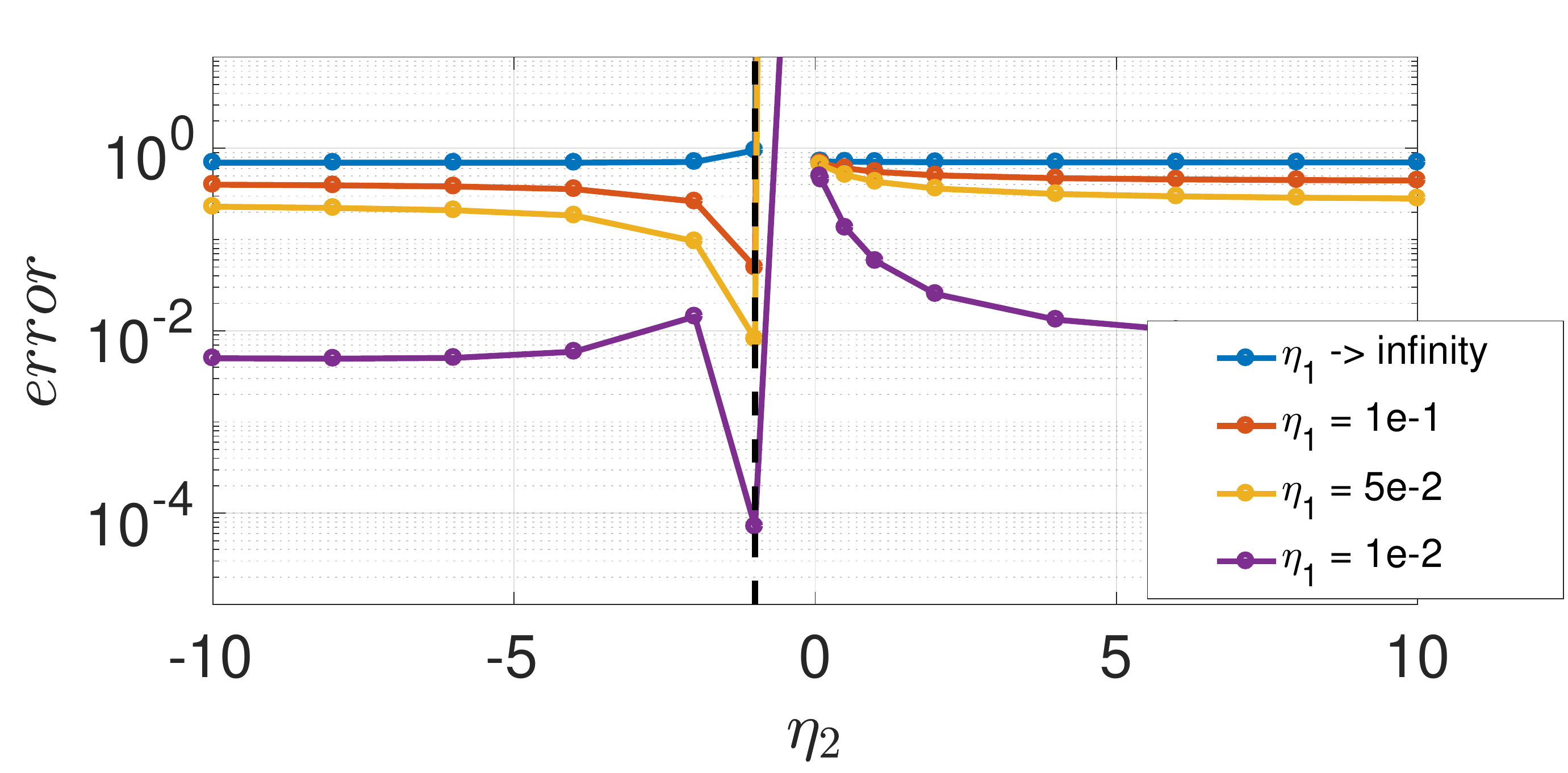}
		\caption{}
	\end{subfigure}
 	\begin{subfigure}{.48\textwidth}
		\centering
 		\includegraphics[width=1\textwidth]{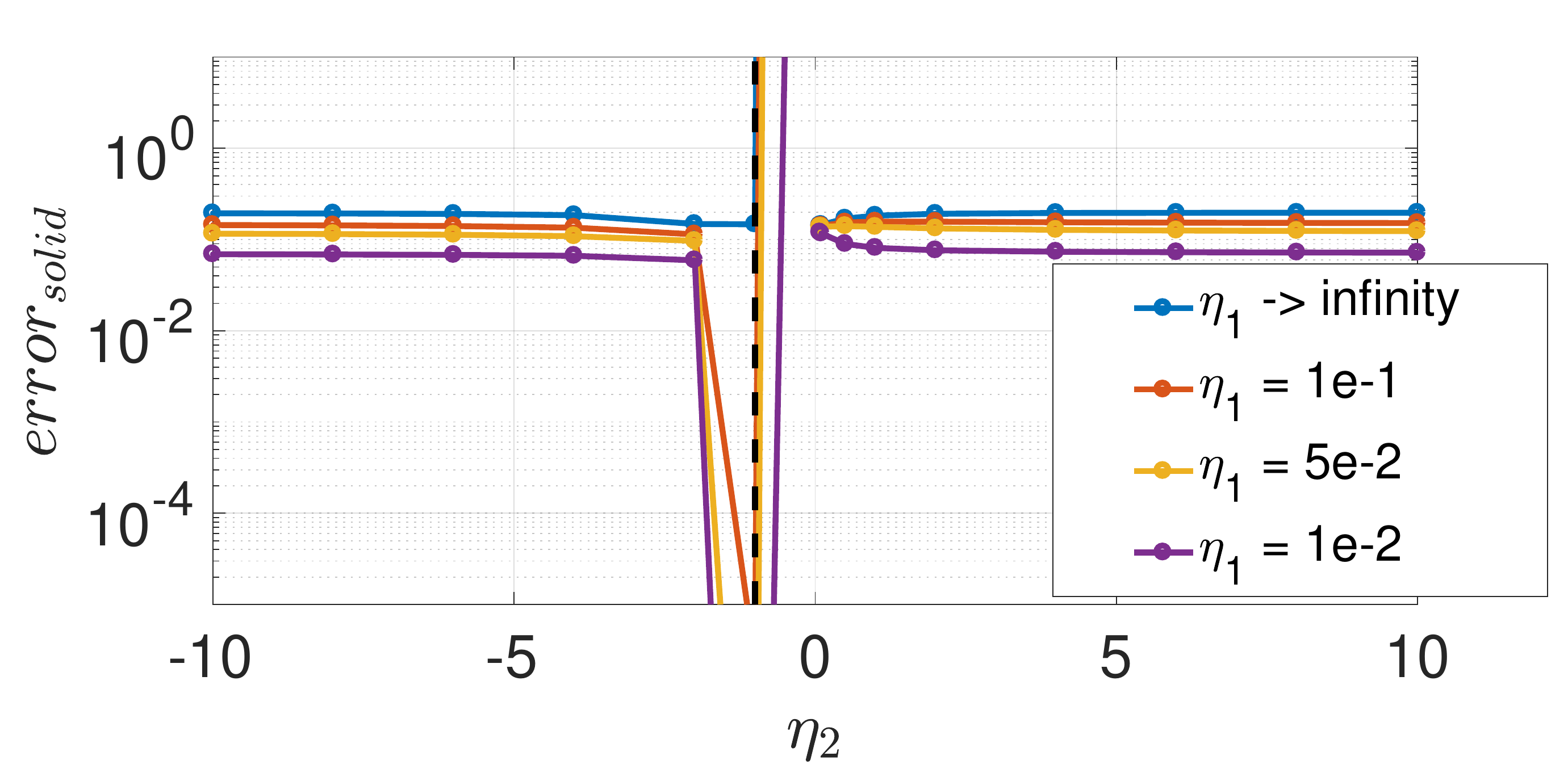}
 		\caption{}
 	\end{subfigure}
	\centering
	\caption{Error comparison for the advection equation, vertical dashed line refers to $\eta_2 = -1/c$, and horizontal dashed line refers to $\eta_2 \rightarrow \infty$. a) Error in the flow ($N=2$). b) Error in the solid, the optimal value is zero ($N=2$). c) Error in the flow ($N=3$). d) Error in the solid, the optimal value is zero ($N=3$). e) Error in the flow (larger penalization parameter, $N=3$). f) Error in the solid (larger penalization parameter, $N=3$).}
	\label{fig:example1}
\end{figure}

Furthermore, to study the effect of $\eta_2$, we run additional simulations for a range of $\eta_2$, and show the errors in Figure \ref{fig:example1}. The errors in the flow and solid regions for $\eta_1 = 10^{-3}$, $\eta_1 = 10^{-4}$, and $\eta_1 = 10^{-5}$ are shown in Figures \ref{fig:example1}a, \ref{fig:example1}b, \ref{fig:example1}c and \ref{fig:example1}d, respectively. For consistency with the analysis of modified equations, the first group of cases is performed in polynomial order $N=2$, and the second group of cases is performed in polynomial order $N=3$. Improved accuracy is seen when the penalization parameter is decreased. In addition, there exists an optimal $\eta_2$ that leads to minimal errors in both the flow and solid regions, which is the same for all penalization parameters. This optimal value is $\eta_2 = -1/c$, indicating that inside the solid the first-order penalization term becomes $-\partial u / \partial x$ thus the physical advection is canceled out. From Figure \ref{fig:example1}b and Figure \ref{fig:example1}d, this cancelation will lead to almost zero error inside the solid, indicating that the boundary condition is satisfied exactly. At a larger $\eta_1$, this optimal value remains valid, but the optimal error increases, as shown in Figure \ref{fig:example1}c and \ref{fig:example1}d. Therefore, to reach the optimal accuracy, we need to use a small penalization parameter $\eta_1$, in combination with the optimal $\eta_2$. These findings are consistent with the theory that the volume penalization modeling error converges with $\eta_1 \rightarrow 0$. Furthermore, the conclusion of the modified equation analysis is also validated, since choosing $\eta_2 = -1/c$ leads to improved accuracy and almost satisfies the boundary conditions exactly.

\begin{figure}[htbp!]
	\begin{subfigure}{.48\textwidth}
		\centering
		\includegraphics[width=165pt]{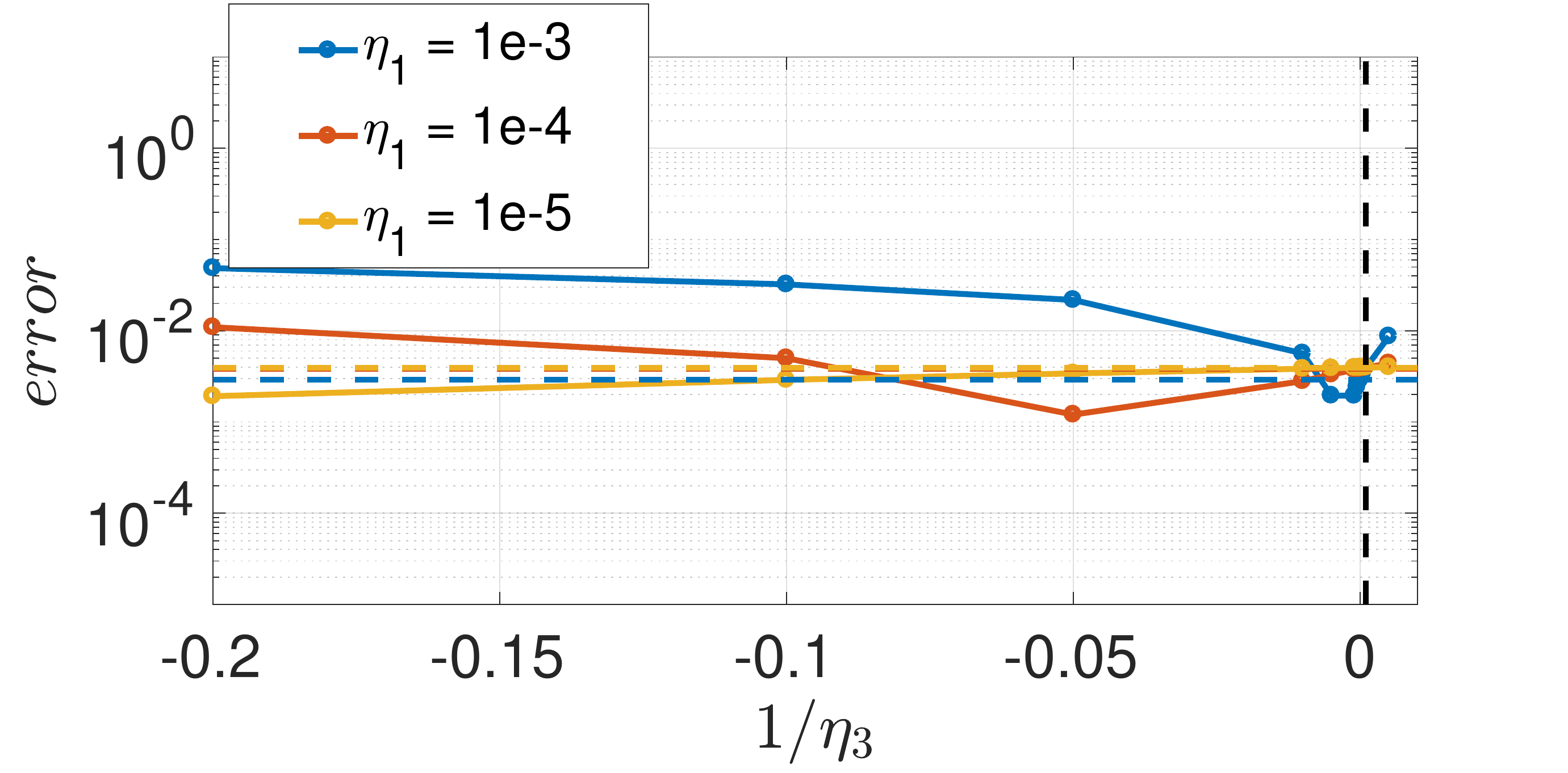}
		\caption{}
	\end{subfigure}
	\begin{subfigure}{.48\textwidth}
		\centering
		\includegraphics[width=165pt]{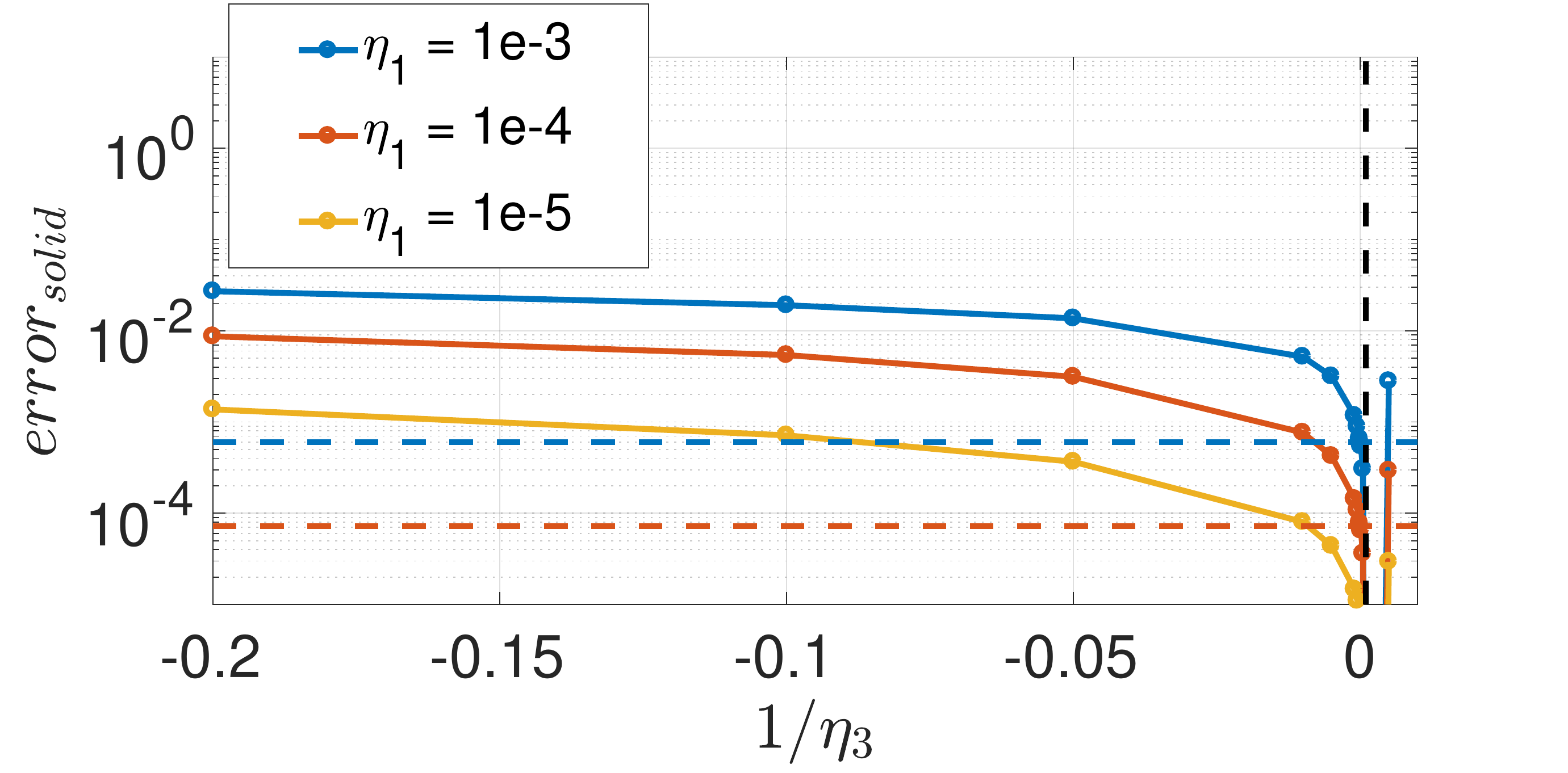}
		\caption{}
	\end{subfigure}
	\begin{subfigure}{.48\textwidth}
		\centering
		\includegraphics[width=165pt]{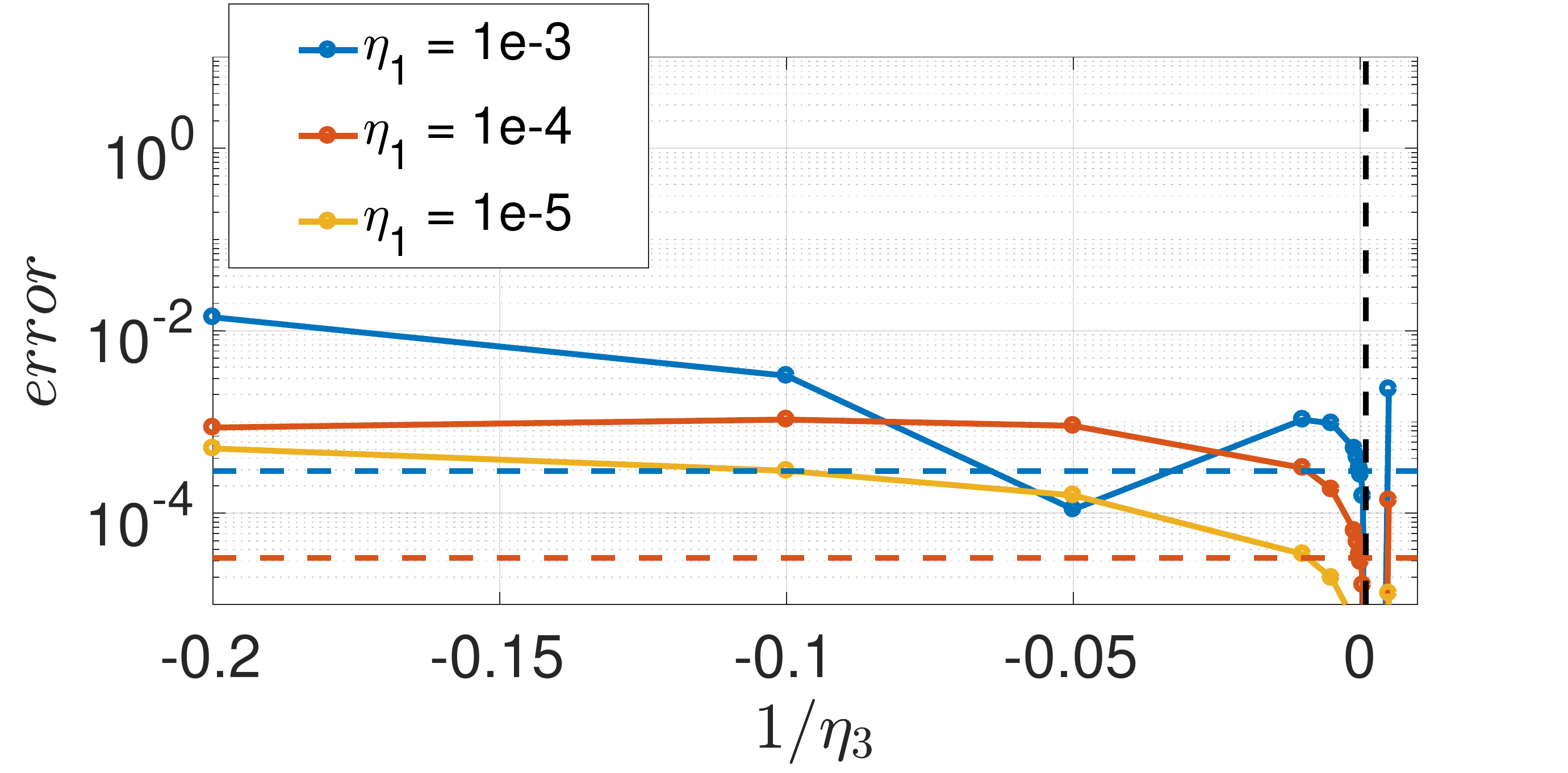}
		\caption{}
	\end{subfigure}
	\begin{subfigure}{.48\textwidth}
		\centering
		\includegraphics[width=165pt]{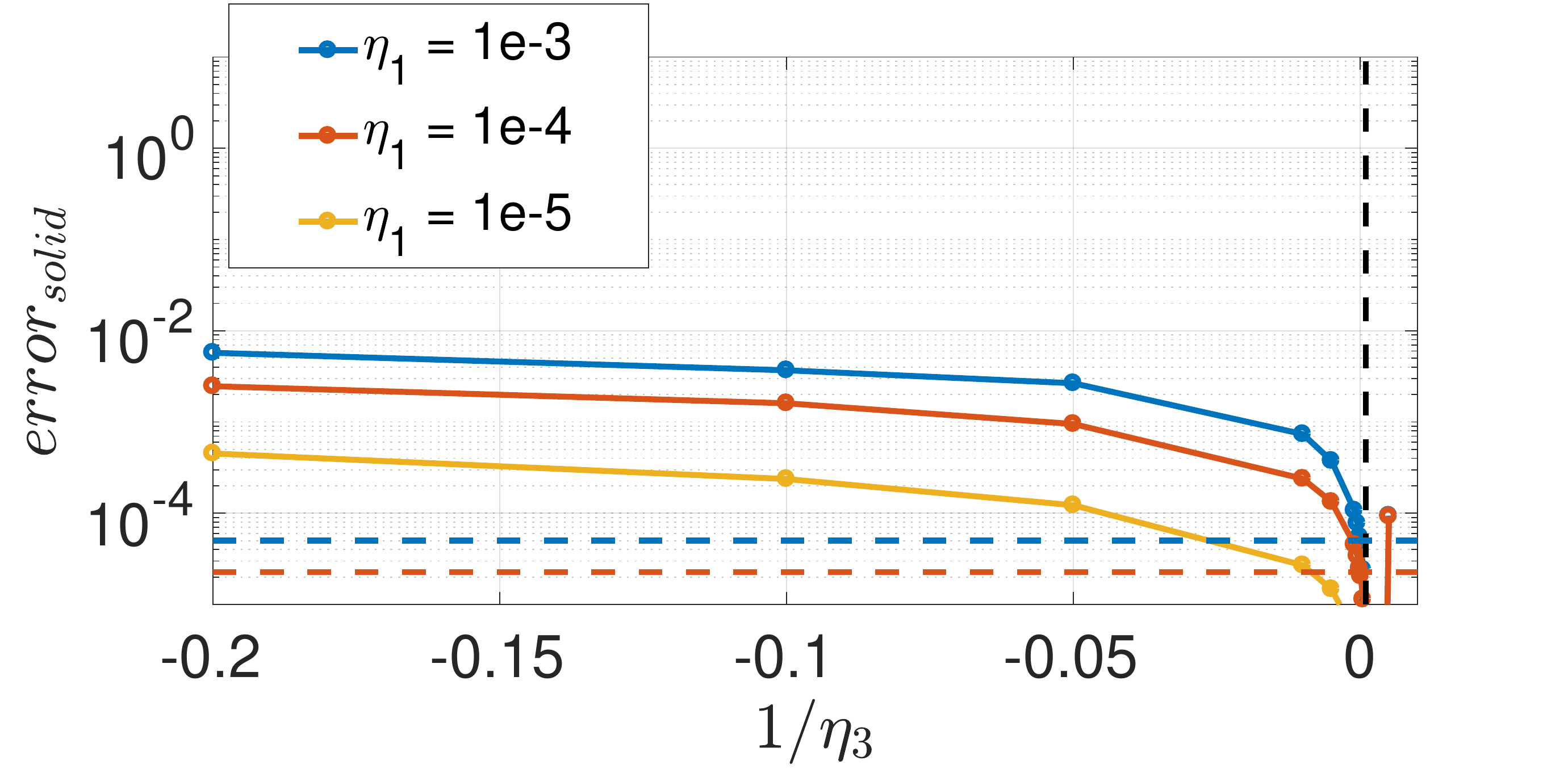}
		\caption{}
	\end{subfigure}
	\centering
	\caption{Error comparison for the advection-diffusion equation ($\nu = 0.001$), vertical dashed line refers to $1 / \eta_3 = \nu = 0.001$, horizontal dashed line refers to $\eta_3 \rightarrow \infty$ (without second-order term) a) Error in the flow ($\eta_2 = -1/c$), BR1. b) Error in the solid ($\eta_2 = -1/c$), BR1. c) Error in the flow ($\eta_2 = -1/c$), LDG. d) Error in the solid ($\eta_2 = -1/c$), LDG.}
	\label{fig:example2}
\end{figure}

\begin{figure}[htbp!]
	\begin{subfigure}{.48\textwidth}
		\centering
		\includegraphics[width=165pt]{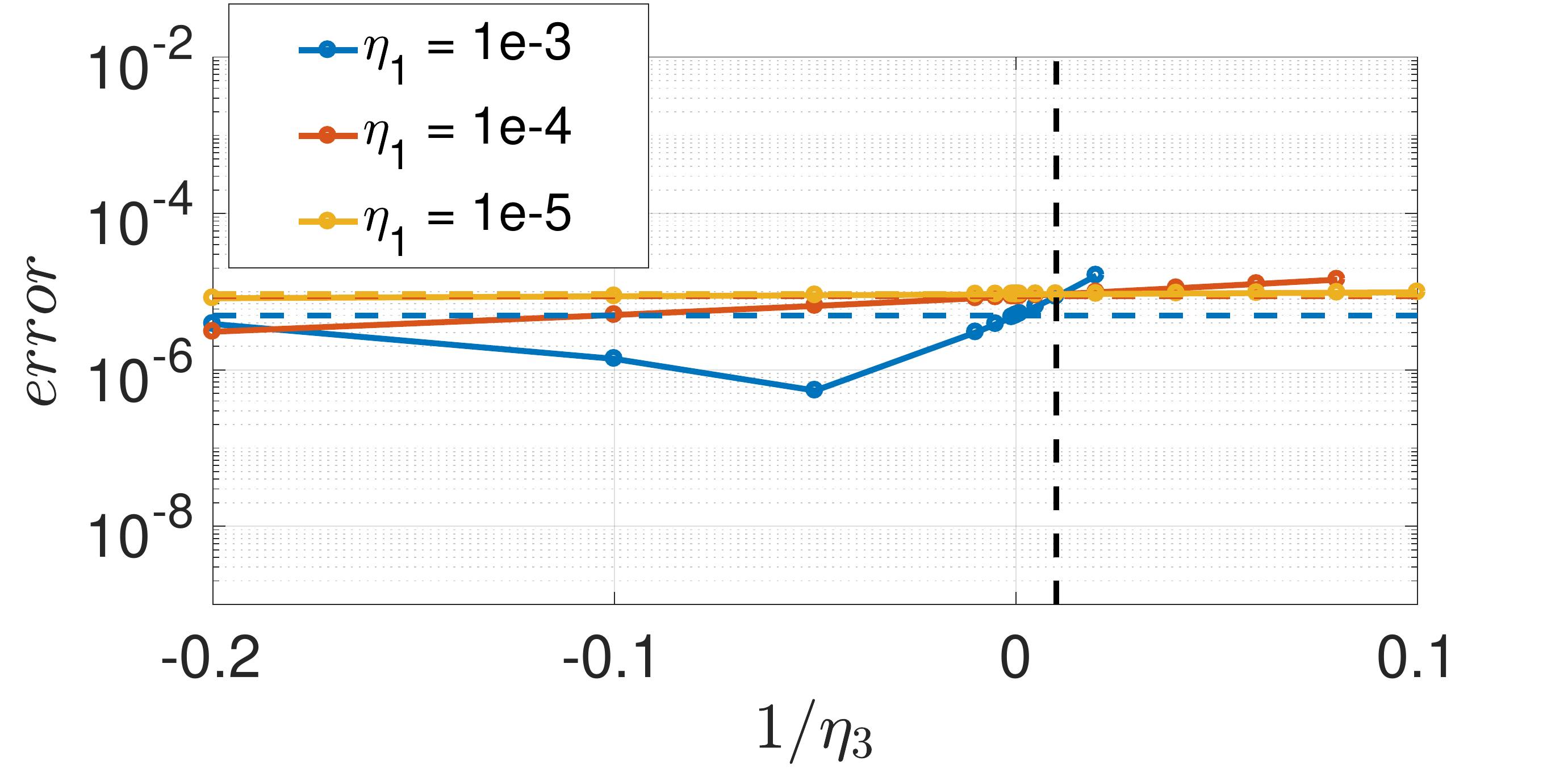}
		\caption{}
	\end{subfigure}
	\begin{subfigure}{.48\textwidth}
		\centering
		\includegraphics[width=165pt]{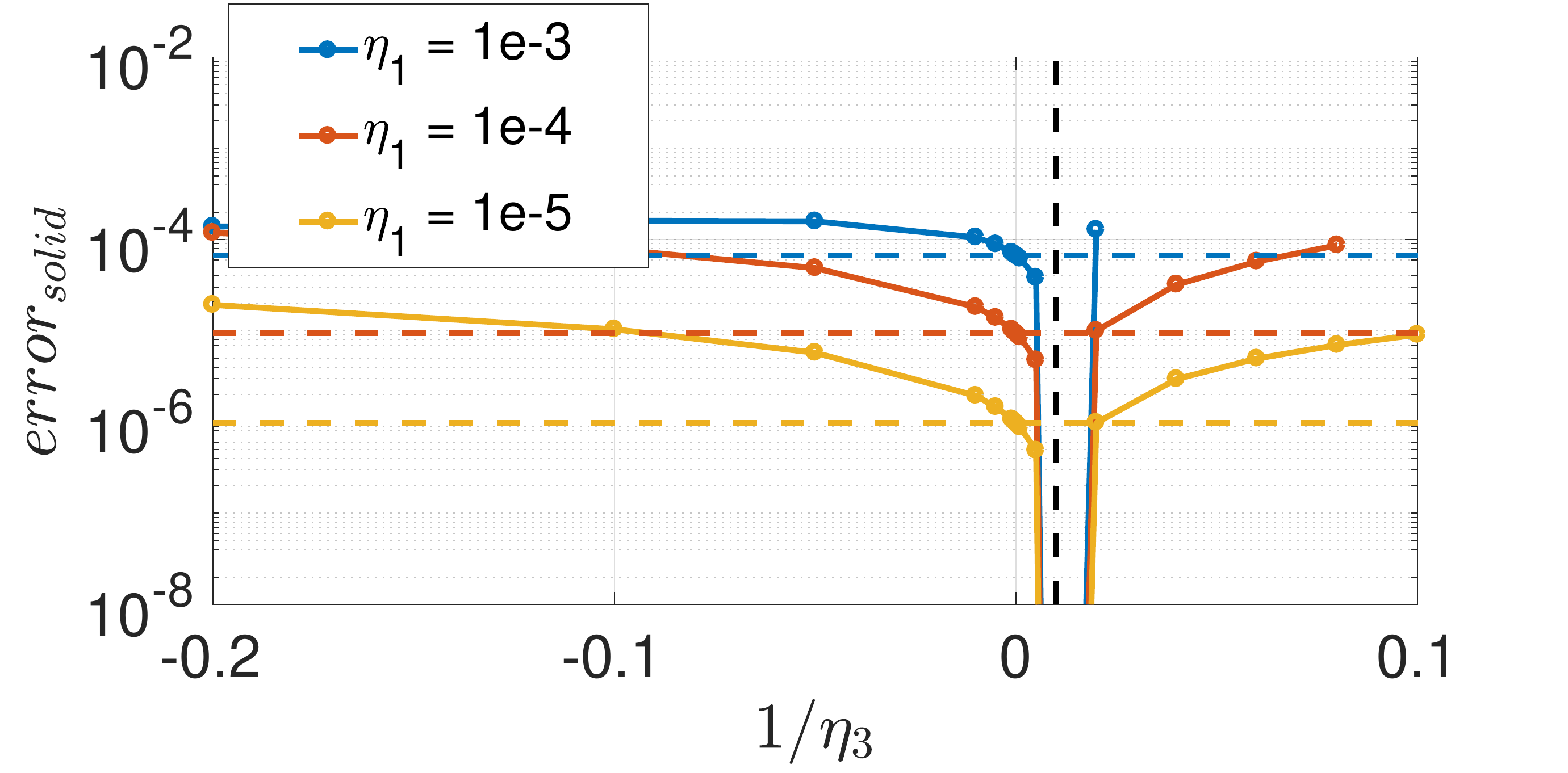}
		\caption{}
	\end{subfigure}
	\begin{subfigure}{.48\textwidth}
		\centering
		\includegraphics[width=165pt]{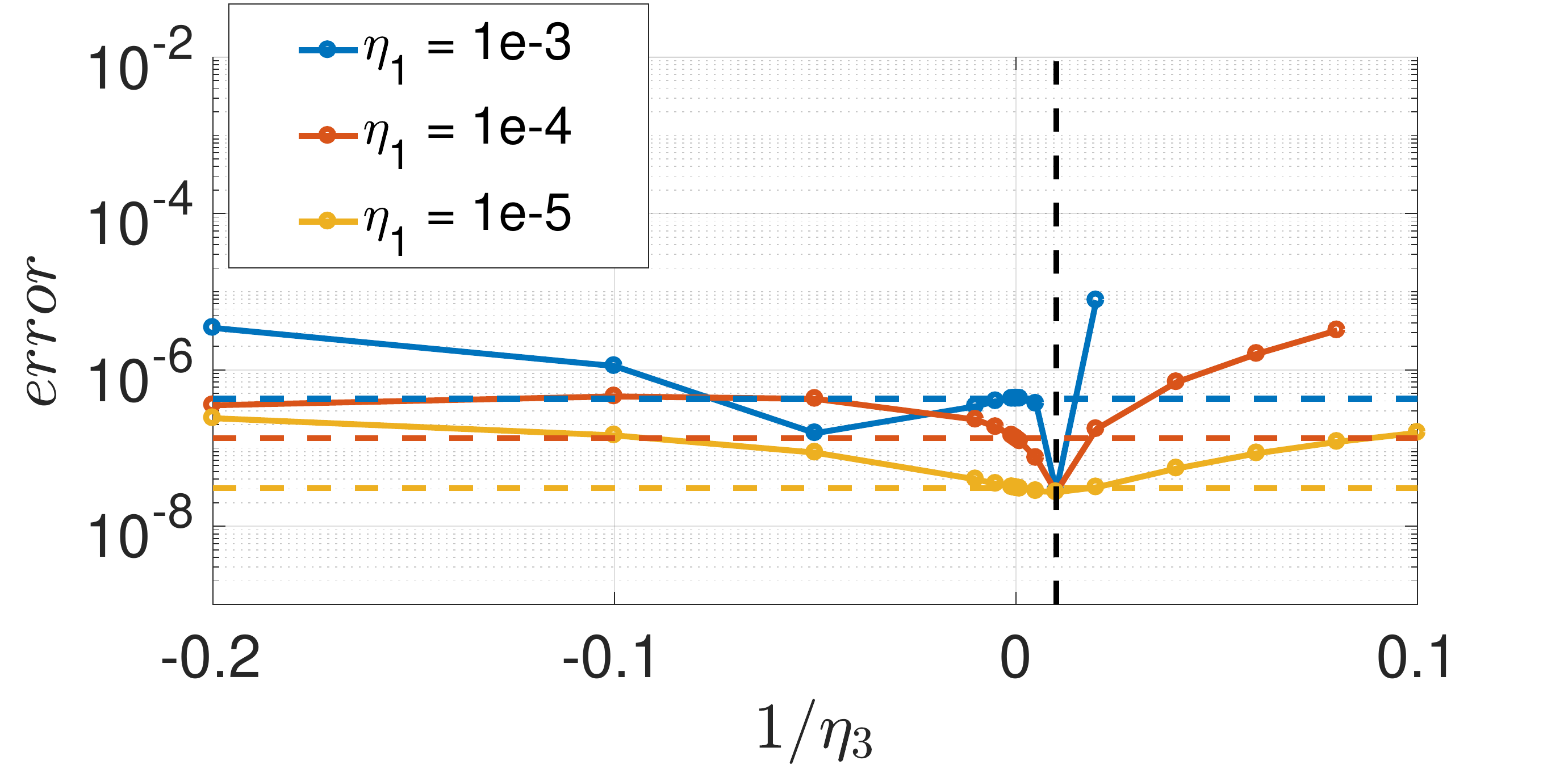}
		\caption{}
	\end{subfigure}
	\begin{subfigure}{.48\textwidth}
		\centering
		\includegraphics[width=165pt]{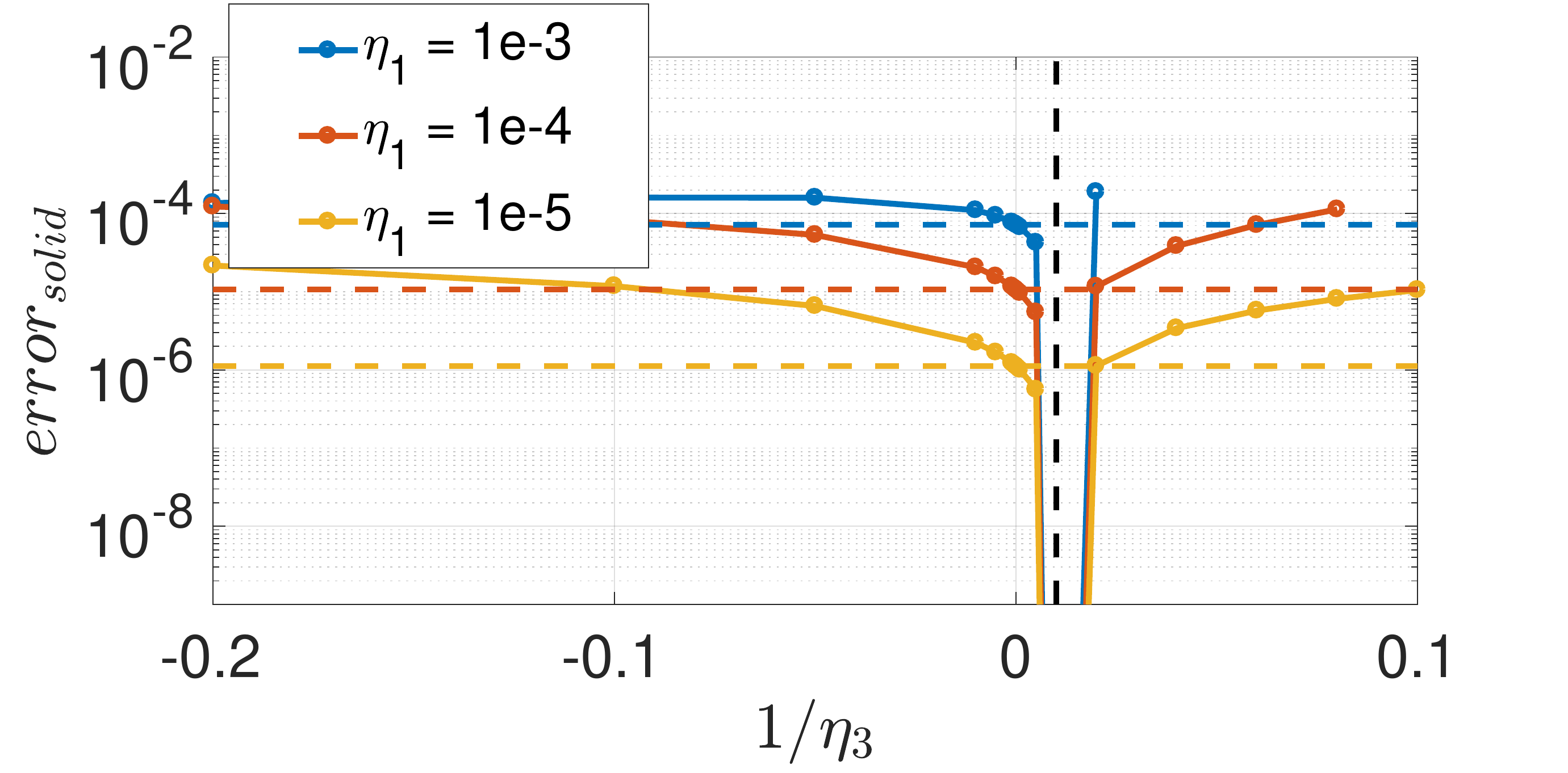}
		\caption{}
	\end{subfigure}
	\centering
	\caption{Error comparison for the advection-diffusion equation ($\nu = 0.01$), vertical dashed line refers to $1 / \eta_3 = \nu = 0.01$, horizontal dashed line refers to $\eta_3 \rightarrow \infty$ (without second-order term). a) Error in the flow ($\eta_2 = -1/c$), BR1. b) Error in the solid ($\eta_2 = -1/c$), BR1. c) Error in the flow ($\eta_2 = -1/c$), LDG. d) Error in the solid ($\eta_2 = -1/c$), LDG.}
	\label{fig:example3}
\end{figure}

To investigate the effect of the viscous term, the advection-diffusion equation is investigated. Since the optimal $\eta_2$ in the advection equation has been obtained, $\eta_2 = -1/c$ is selected for all cases. We proceed as for the advection equation, we set $K = 40$ elements, $\Delta_{\text{s}} = \Delta x$, $r = 1/40$ and $N = 3$. The initial condition with wavenumber $\overline{\omega} \Delta x/(N+1) = 0.3223$ is used and marched in time to $t = 1.5$. Taking into account the effect of diffusion, the error in the flow region is limited to $x \in [\Delta_{\text{s}}, 0.7]$. 

For the discretization of the viscous flux, either the BR1 or the LDG scheme is considered. The results for two physical viscosities $\nu = 0.001$ and $\nu = 0.01$ are shown in Figures \ref{fig:example2} and Figures \ref{fig:example3}, respectively. For both cases, it is observed that the optimal second-order coefficient $\eta_3$ exists and can lead to a minimal error within the solid (as shown in Figure \ref{fig:example2}b and \ref{fig:example2}d and Figure \ref{fig:example3}b and \ref{fig:example3}d). This optimal value shows the relationship $1 / \eta_3 = \nu$, which also indicates the cancelation of the viscous term inside the solid. This agrees with the optimal $\eta_3$ derived from the modified equation analysis. However, when looking at the error inside the flow, the optimal second-order penalization term does not lead to the lowest error when the BR1 scheme is used. This highlights the importance of choosing appropriate Riemann solvers to maintain good accuracy in the flow region. When the LDG scheme is selected, the optimal $\eta_3$ will reach the lowest error in the flow region, indicating that this flux is more suitable for the present problem. Therefore, when handling the viscous term, the LDG scheme is preferred, which gives consistent results of errors in the solid and in the fluid. In summary, the one-dimensional test case shows that the optimal penalization parameters derived from the modified equation analysis achieve minimal numerical errors in imposing the boundary conditions.

\subsection{Two-dimensional advection-diffusion equation}\label{subsec:2DTest}
\begin{figure*}[htbp]
		\centering
		\includegraphics[width=150pt]{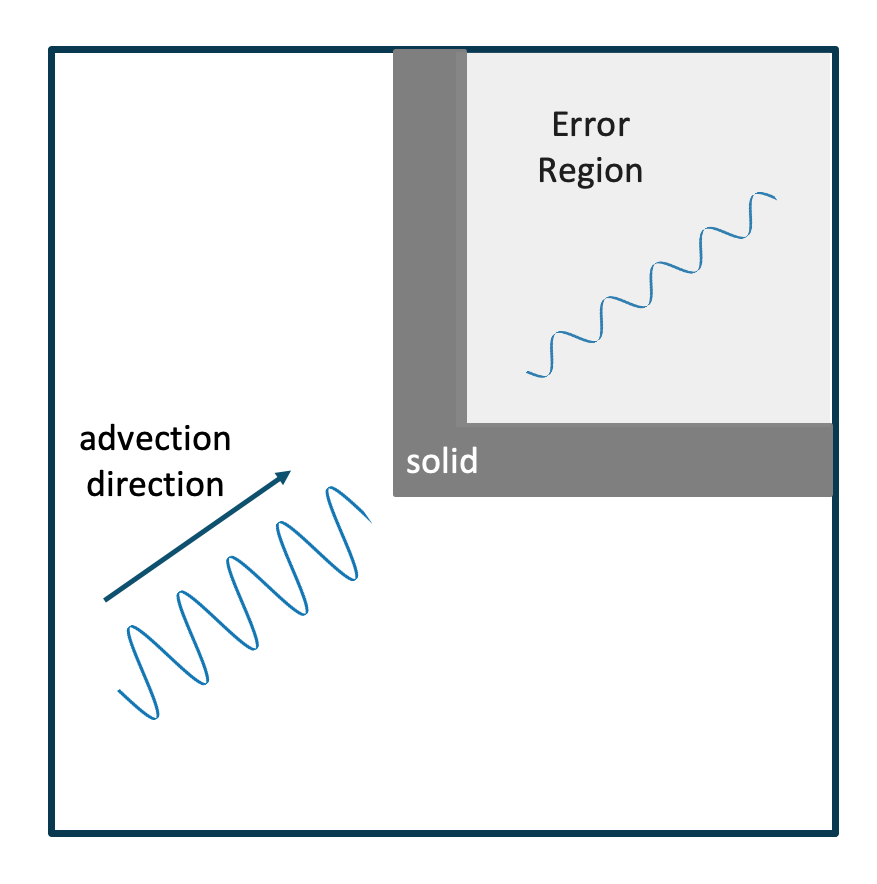}
		\caption{Schematic illustration of the advection problem with IBM.}
		\label{fig:cartoon2}
\end{figure*}

In this section, a numerical experiment is performed for the two-dimensional advection-diffusion equation, using the conclusions of the modified equation analysis. The one-dimensional test case in the previous section is extended to two space directions. Therefore, the optimal parameters derived from one-dimensional test cases are then dependent on each space direction.  Extensions for GL points can be found in \cite{hesthaven2007nodal}. The governing equation is (the solid region is the non-slip wall):

\begin{align}
\frac{\partial u}{\partial t} + \nabla \cdot (\mathbf{f}_{\text{adv}} + \mathbf{f}_{\text{diff}}) + \frac{\chi}{\eta_1} u + \nabla\cdot\left(\mathbf{g}\chi u\right) + \nabla\cdot\left(\mathbf{H}\nabla\left(\chi u\right)\right) = 0,
\end{align}
where the advection flux is $\mathbf{f}_{\text{adv}} = (c_x u,c_y u)^T$, the diffusion flux $\mathbf{f}_{\text{diff}} = (-\nu_x\partial u/\partial x,-\nu_y\partial u/\partial y)^T$, $\mathbf{g} = (1/\eta_{2,x},1/\eta_{2,y})^T$, and $\mathbf{H} = \text{diag}\left(1/\eta_{3,x},1/\eta_{3,y}\right)$. The first-order and second-order penalization parameters in each direction is denoted by the second subscript. Here we set $c_x = c_y = 1$ and $\nu_x = \nu_y = 0.001 $, therefore, the optimal parameters satisfy $\eta_{2,x} = \eta_{2,y}$ and $\eta_{3,x} = \eta_{3,y}$. Note that, for the present linear equation, the extension to different advection velocities and viscosities in each direction is straightforward, while the optimal penalization parameter (i.e., trivial solution from modified equation analysis) also varies in different directions. As in the one-dimensional test case, the solid wall is considered in the middle of the computational domain. A schematic illustration of the two-dimensional problem for the present study is shown in Figure \ref{fig:cartoon2}. The solid non-slip region has a $L$ shape, which is centered in the middle of the square domain, making the top right region amplified by the wall. If the advection direction is set appropriately, the initial wave will move towards the wall. After that, we can solve the equation until all the solutions in the top right region have been penalized (which are then expected to be zero), and compute the error in this region. The error is again the difference between the numerical and exact solution (here set to zero). 

We consider a square computational domain in $x \in [-0.1,0.1]$ and $y \in [-0.1,0.1]$, with periodic boundary conditions. The domain is discretized into $20$ equispaced elements in both the $x$ and the $y$ directions, resulting in $400$ square elements in total and uniform mesh size $\Delta x = \Delta y = 0.01$. The penalization parameter and the explicit time step is set to $\eta_1 = \Delta t = 10^{-4}$. The polynomial order $N=3$ is selected. Due to the preset flow advection parameters, the advection moves towards the top right direction. The width of the solid region is set to the size of a uniform grid $\Delta_{\text{s}} = \Delta x$, resulting in the solid ratio $r = 1/20$. We use the wavelike initial condition $u(x,y) = \sin(\omega x + \omega y)$, where a nondimensional wavelength $\overline{\omega} \Delta x/(N+1) = 0.3307$ is selected. 

\begin{figure*}[htbp]
    \begin{subfigure}{.32\textwidth}
 		\centering
		\includegraphics[width=115pt]{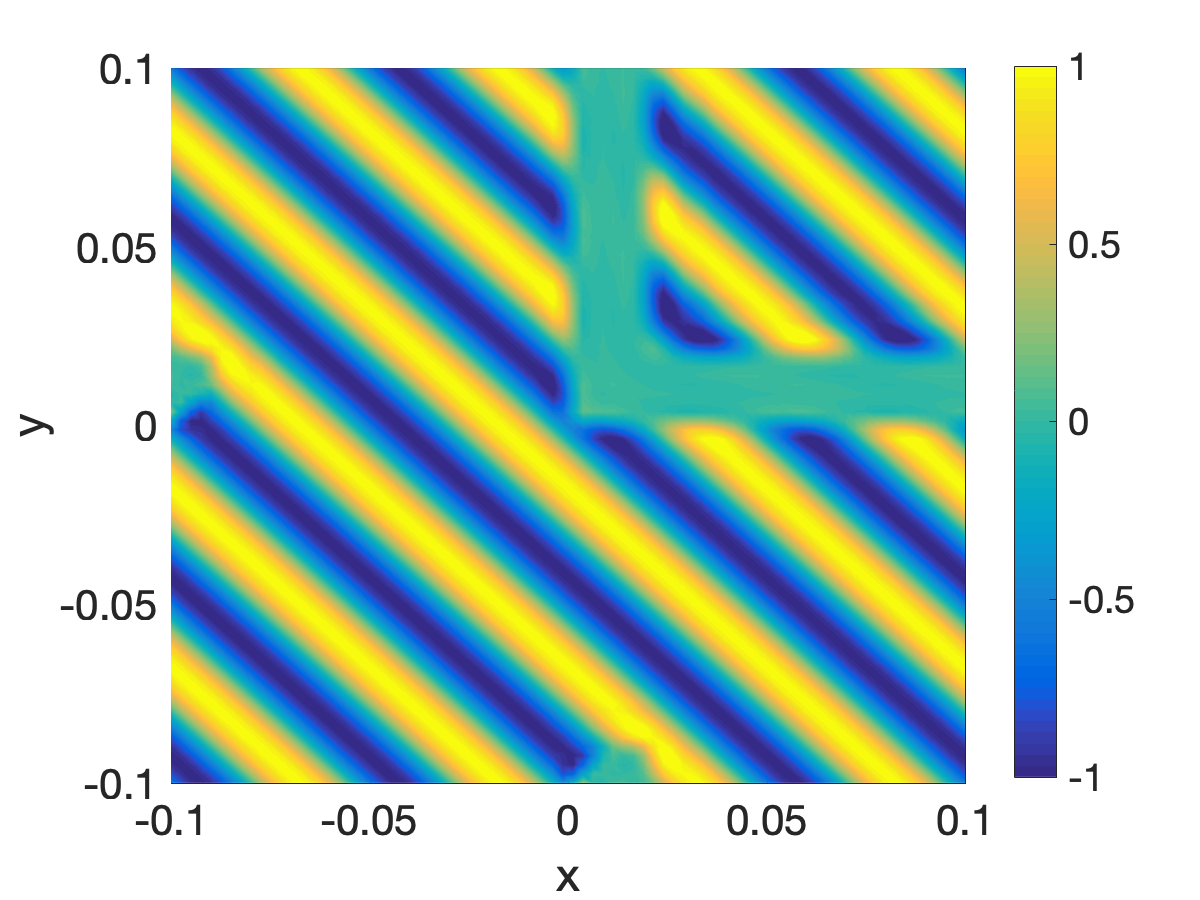}
		\caption{}
	\end{subfigure}
	\begin{subfigure}{.32\textwidth}
 		\centering
		\includegraphics[width=115pt]{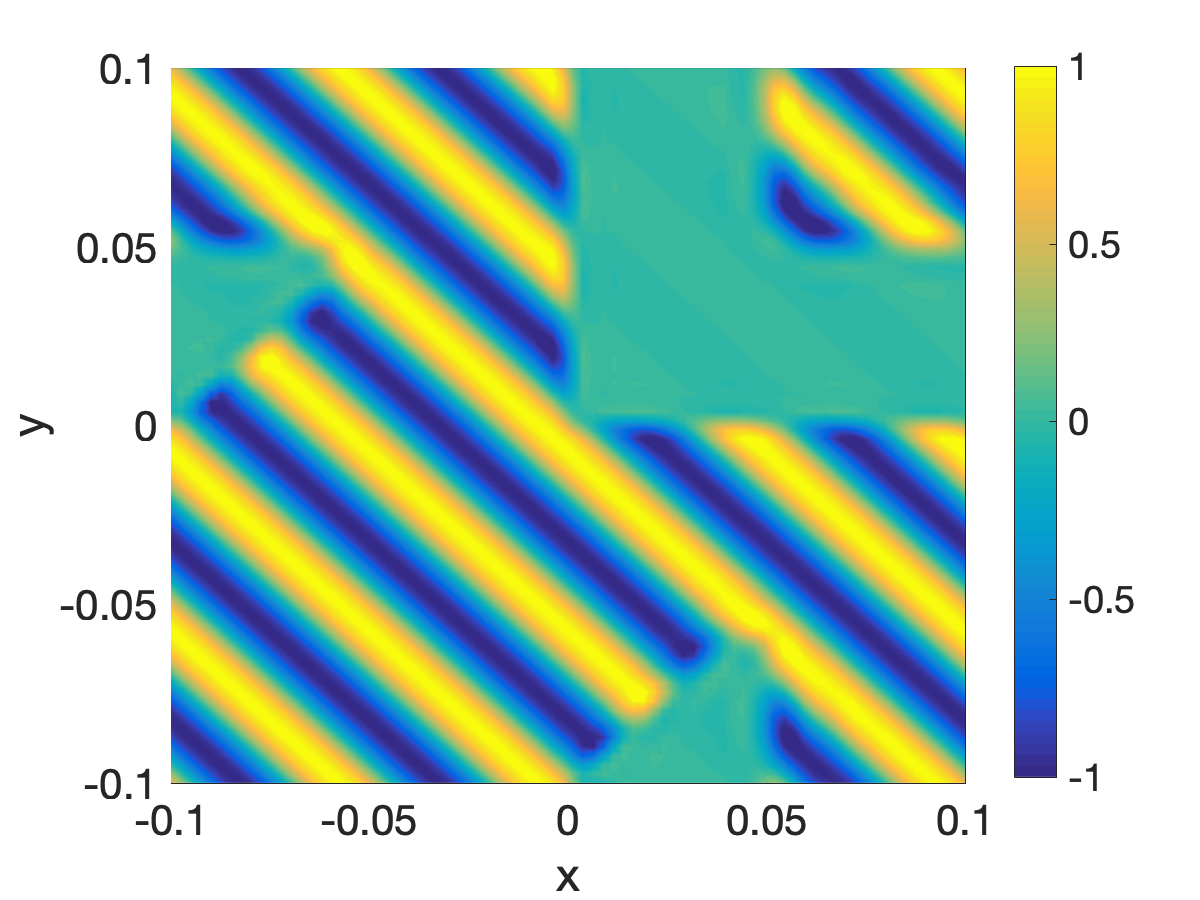}
		\caption{}
	\end{subfigure}
	\begin{subfigure}{.32\textwidth}
 		\centering
		\includegraphics[width=115pt]{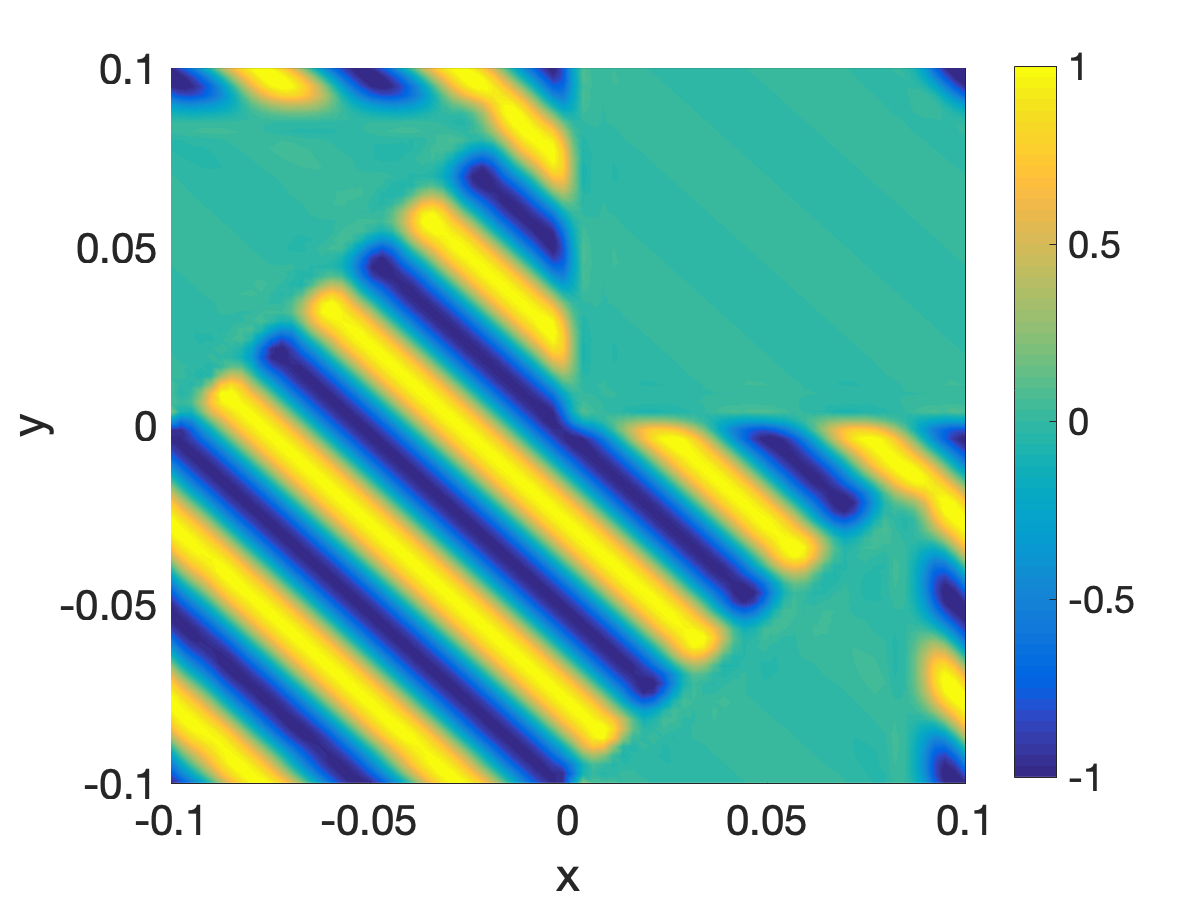}
		\caption{}
	\end{subfigure}
	\centering
	\caption{Simulation under different parameters ($r = 1/20$, $N = 3$, initial wavenumber $\overline{\omega} \Delta x/(N+1) = 0.3307$, $K = 20$). a) t = 0.01. b) t = 0.04. c) t = 0.08.}
	\label{fig:advection-field}
\end{figure*}
\begin{figure*}[htbp]
    \begin{subfigure}{.48\textwidth}
 		\centering
		\includegraphics[width=165pt]{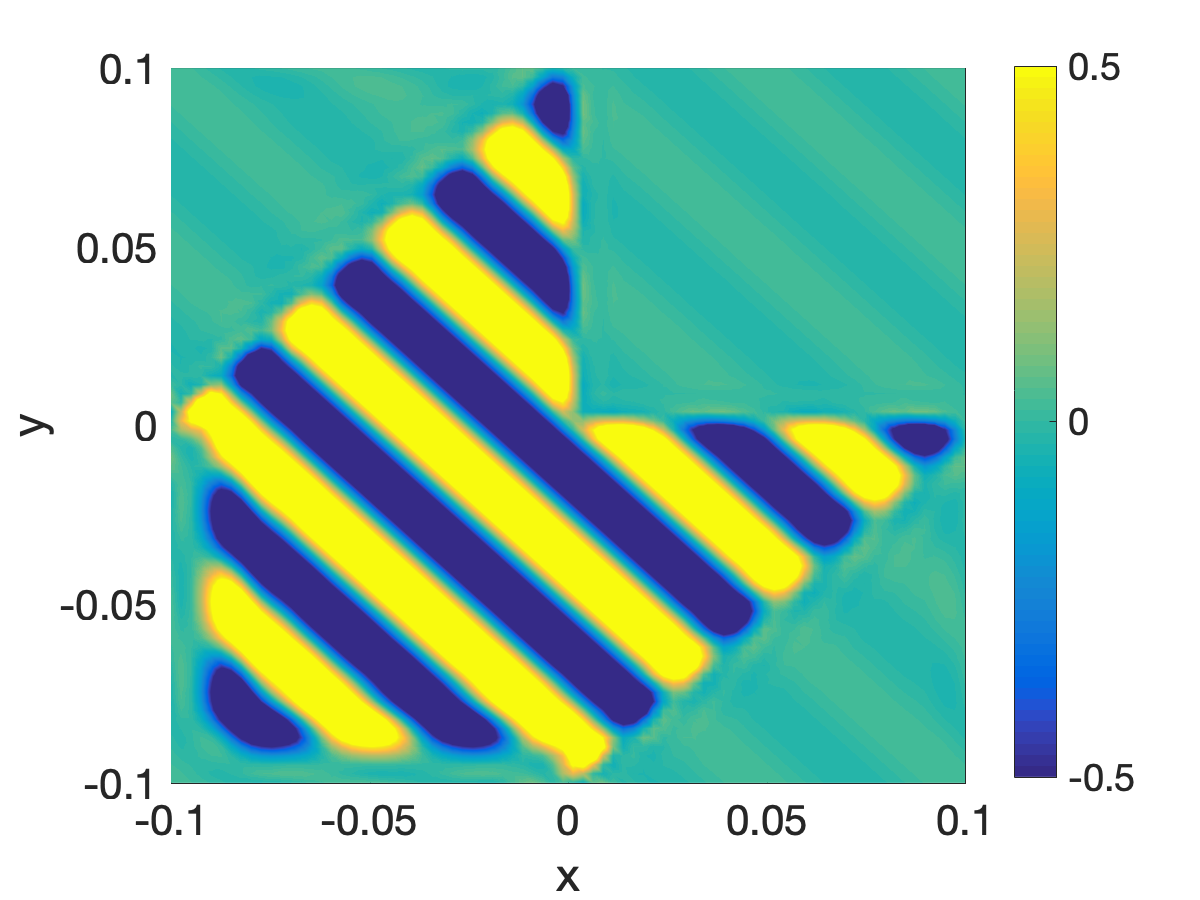}
		\caption{}
	\end{subfigure}
	\begin{subfigure}{.48\textwidth}
 		\centering
		\includegraphics[width=165pt]{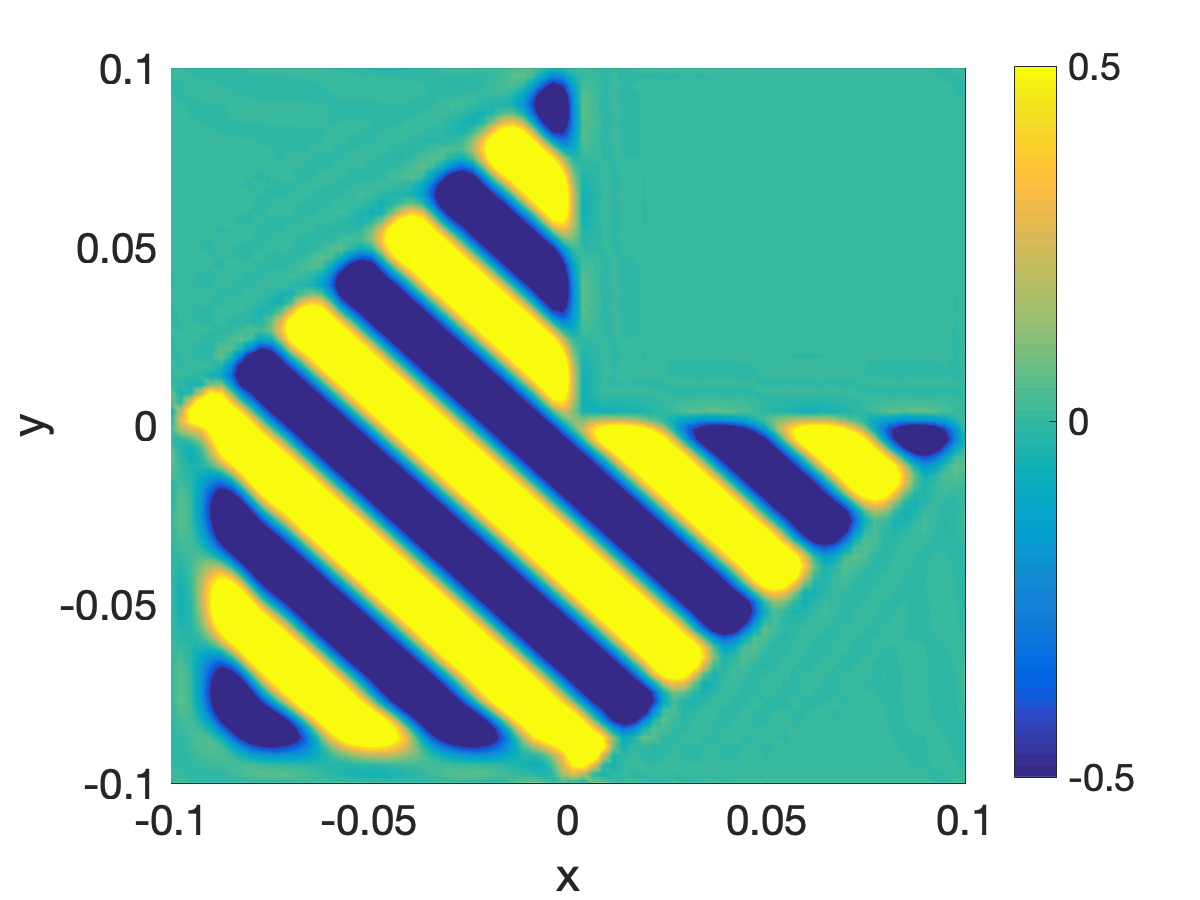}
		\caption{}
	\end{subfigure}
	\centering
	\caption{Simulation under different parameters ($r = 1/20$, $N = 3$, initial wavenumber $\overline{\omega} \Delta x/(N+1) = 0.3307$, $K = 20$). The difference lies in the upper right flow region. a) $\text{error} = 0.0207$, $\text{error}_{\text{solid}} = 0.0552$. b) $\text{error} = 1.4616\cdot 10^{-5}$, $\text{error}_{\text{solid}} = 0$. }
	\label{fig:advection-field2}
\end{figure*}

The first simulation for only advection is performed when only the first penalization term ($\eta_1$ for $u$) is included. Figure \ref{fig:advection-field} shows three typical solution fields at different times. As shown in the figure, the penalized solution will move towards the top right corner, and finally dominate the entire domain due to the periodic boundary conditions. To compare the accuracy of simulation, the final simulation time is set to $0.11$. Two solution fields, without and with the optimal first-order penalization term, are shown in Figure \ref{fig:advection-field2}, where the values of $\eta_2,x$ and $\eta_2,y$ are set to $-1$ to match the physical advection speed. The improved accuracy from adding the optimal first-order term is seen in both the solution field and in the error. The error in the fluid region has been greatly reduced from $0.0207$ to $1.4616\cdot 10^{-5}$. This numerical experiment extends and validates the conclusions obtained from the modified equation analysis, where the optimal first-order penalization term cancels the advection term and leads to improved accuracy.

\begin{table}[htbp!]
	\vspace{20pt}
	\centering
	\begin{tabular}{cccc}%p{2.5cm}p{2.5cm}p{2.5cm}p{4cm}
		\hline
		Diffusive Flux Scheme & $\eta_1 = 10^{-4}$ & \makecell[l]{$\eta_1 = 10^{-4}$ \\ $\eta_2 = -1$} & \makecell[l]{$\eta_1 = 10^{-4}$ \\ $\eta_2 = -1$ \\ $\eta_3 = 10^{3}$} \\
		\hline
		BR1 & $1.6610\cdot 10^{-4}$ & $1.3513\cdot 10^{-4}$ & $1.5874\cdot 10^{-4}$ \\
		LDG & $6.4091\cdot 10^{-5}$ & $7.1993\cdot 10^{-6}$ & $2.2669\cdot 10^{-7}$ \\  
		\hline       
	\end{tabular}
	\caption{Error comparison of the advection equation with IBM wall under different diffusive flux schemes and different combinations of penalization parameters.}
	\label{table2}
\end{table}

Additional numerical experiments are performed for the advection-diffusion equation. The space and time discretizations remain the same as in the advection case. The final time is set to $t = 0.15$. Three strategies are considered: 1) only volume penalization for the non-slip wall boundary condition, 2) volume penalization for both the value and the first-order term, and 3) volume penalization for all terms. Two types of viscous fluxes with either the BR1 or the LDG scheme are considered. The errors inside the fluid region are compared in Table \ref{table2}, where conclusions similar to one-dimensional advection can be drawn. When the BR1 scheme is used, adding additional first-order and second-order penalization terms improves the overall accuracy, compared with the standard case (the first strategy). However, the addition of a second-order term does not lead to improved accuracy in the flow region. When the LDG scheme is used, adding first- and second-order terms will lead to a greater reduction of the error. This is consistent with the observations for the one-dimensional advection equation, where the LDG scheme is shown to provide more accurate results than the BR1 scheme. This numerical experiment validates the proposed modified equation analysis for the second-order derivative in two-dimensional linear equations.

\section{Conclusions}\label{sec:Conclusion}
This study contributes to a better understanding of the numerical errors for Immersed Boundary Methods based on volume penalization, in combination with a high-order nodal discontinuous Galerkin scheme. For this purpose, an analysis of the modified equation is provided.

The modified equation is a useful tool to analyze dissipative/dispersive errors related to the numerical discretizations. In this paper, we focus on the spatial errors introduced by the Immersed Boundary Method. Nodal solutions are expanded as Taylor series, and by rearranging the pseudo-differential equation new terms arise. These terms allow us to obtain insight into the dissipative/dispersive character of the errors and guidelines for their minimization. For example, the inclusion of extra penalization terms of the first and second derivatives, in addition to the classic penalization of the variable. Through this analysis, we provide optimal values for the first- and second-order penalization parameters to cancel the advection/diffusive errors inside the solid, which in turn lead to improved errors in the flow.

Numerical experiments validate the theoretical findings obtained from the analysis of modified equations, where optimal penalization parameters can lead to minimal errors (with a sufficiently small penalization parameter $\eta_1$). When combined with an appropriate numerical scheme (here, Local discontinuous Galerkin for viscous terms), minimal errors in the flow region are reached.

Future work will extend these findings to systems of partial differential equations with non-linearities, and extend the theoretical analysis to multi-dimensional systems.

%This paper presents a modified equation analysis for Immersed Boundary Methods based on volume penalization, in combination with a high-order nodal discontinuous Galerkin scheme. We analyse the performance of volume penalization for a linear advection-diffusion equation. Several conclusions can be drawn: CAMBIAR ESTO
%\begin{enumerate}
%    \item The modified equation analysis is a useful tool to analyse the error of IBM treatment. It can be performed based on different IBM approaches and different numerical schemes, to obtain insights and guidelines for the implementation of immersed boundary treatment.
%    \item Through theoretical analysis, we provide a family of solutions that lead to improvements in IBM treatment. In particular, a trivial solution is obtained for optimal first-order and second-order penalization parameters, which cancels the advection and diffusive errors inside the solid.
%    \item Numerical experiments validate the conclusions from the modified equation analysis, where the optimal penalization parameters can lead to minimal errors inside the solid (with a sufficiently small penalization parameter $\eta_1$). When combined with an appropriate numerical scheme (here LDG), minimal errors in the flow region can be reached.
%\end{enumerate}
%Future work will extend these findings to a system of PDEs and nonlinear equations, as well as investigate non-trivial solutions for different DG treatments. 

\begin{acknowledgements}
JK and EF acknowledge the financial support of the European Union’s Horizon 2020 research and innovation programme under the Marie Skłodowska-Curie grant agreement (MSCA ITN-EID-GA ASIMIA No. 813605). VJL, EF, and EV acknowledge financial support from the European High-Performance Computing Joint Undertaking (JU) under grant agreement (No. 956104). The JU receives support from the European Union’s Horizon 2020 research and innovation programme under grant agreement (No. 823844) and Spain, France, Germany.
EF would like to thank the support of the Spanish Ministry MCIN/AEI/10.13039/501100011033 and the European Union NextGenerationEU/PRTR for the grant ``Europa Investigación 2020'' EIN2020-112255, and also the Comunidad de Madrid through the call Research Grants for Young Investigators from the Universidad Politécnica de Madrid. 
Finally, all authors gratefully acknowledge the Universidad Politécnica de Madrid (www.upm.es) for providing computing resources on Magerit Supercomputer.
\end{acknowledgements}

\section*{Declarations}
% \textbf{Declarations}

\textbf{Conflict of interest} The authors declare no competing interests.
% \section*{Code Availability}

\noindent
\textbf{Code Availability}
The codes developed during the current study are available from the corresponding author on request.

%There is no conflict of interest.

\appendix

\section{The DGSEM technique}\label{appx:DGSEM}
We re-write Equation (\ref{eqn:PDE}) in its weak form:

\begin{align}
\int_{0}^{L}\left(\frac{\partial u}{\partial t} + \frac{\partial \widehat{f}}{\partial x} + \frac{\chi}{\eta_{1}} u\right)\psi\,\diff x = 0,
\end{align}
where $\psi = \psi(x,t)$ is a local smooth test function. Given that $\Omega = [0,L]$ is divided into $K$ elements, the integral is split into the sum of element integrals:

\begin{align}
\sum_{k=1}^{K}\left\{\int_{x_{k-1}}^{x_{k}}\left(\frac{\partial u}{\partial t} + \frac{\partial \widehat{f}}{\partial x} + \frac{\chi}{\eta_{1}} u\right)\psi\,\diff x\right\} = 0.
\end{align}
Each element, $x = x(\xi)$, is transformed according to: $x = x_{k-1} + (\xi+1)\Delta x_k/2$, where $\Delta x_k = x_{k} - x_{k-1}$ and $-1\leq\xi\leq 1$. Then, $\diff x = (\Delta x_k/2)\diff\xi$ and $\partial/\partial x = (2/\Delta x_k)\partial/\partial\xi$. Thus the weak form becomes:

\begin{align}
\sum_{k=1}^{K}\left\{\frac{\Delta x_k}{2}\int_{-1}^{1}\left(\frac{\partial u}{\partial t} + \frac{2}{\Delta x_k}\frac{\partial \widehat{f}}{\partial \xi} + \frac{\chi}{\eta_{1}} u\right)\psi\,\diff\xi\right\} = 0.
\end{align}
Assuming that global variables are represented by $K$ local polynomial variables and substituting the Lagrange interpolation of the test function into the Galerkin weak form, $\psi = \sum\psi_jl_j$, we get the following:

\begin{align}
\frac{\Delta x_{k}}{2} \int_{-1}^{1} l_{j}\frac{\partial u_{h}^{k}}{\partial t} \diff \xi + \int_{-1}^{1} l_j\frac{\partial\widehat{f}_{h}^{k}}{\partial\xi} \diff \xi + \frac{\Delta x_{k}}{2\eta_1} \int_{-1}^{1} l_j\chi^ku_h^k \diff \xi = 0,
\end{align}
for $k = 1,2,\ldots,K$ and $j = 0,1,\ldots,N$. The first and third integrals are evaluated as follows:

\begin{align}
& \int_{-1}^{1} l_{j}\frac{\partial u_{h}^{k}}{\partial t} \diff \xi = \sum_{i=0}^{N}\int_{-1}^{1} l_{i}l_{j} \diff \xi \frac{\diff u_{h,i}^{k}}{\diff t}, \label{eqn:integral1} \\
& \int_{-1}^{1} l_j\chi^ku_h^k \diff \xi = \sum_{i=0}^{N}\int_{-1}^{1} \chi^{k}l_{i}l_{j} \diff \xi u_{h,i}^{k}, \label{eqn:integral2}
\end{align} 
whereas the second integral is integrated by parts,

\begin{align}
\int_{-1}^{1} l_j\frac{\partial\widehat{f}_{h}^{k}}{\partial\xi} \diff \xi = \left.l_{j}\flux^{k}\right|_{-1}^{1} - \int_{-1}^{1} l'_j\widehat{f}_{h}^{k} \diff \xi, \label{eqn:integral3}
\end{align} 
being $l'_j = \diff l_j/ \diff \xi$. The VP flux function in the first term is substituted by a numerical flux, i.e.,

\begin{align}
& \flux_1^k \coloneqq \flux(u_{h,N}^k,u_{h,0}^{k+1};+\mathbf{e}^k_x), \\
& \flux_{-1}^k \coloneqq \flux(u_{h,N}^{k-1},u_{h,0}^k;-\mathbf{e}^k_x),
\end{align}
depending on the normal at the boundary, $\pm \mathbf{e}^k_x$, and the solution at two adjacent elements. We discuss the choice of the numerical flux later. The remaining integral is divided as follows:

\begin{align}
\int_{-1}^{1} l'_j\widehat{f}_{h}^{k} \diff \xi = \sum_{i=0}^{N}\int_{-1}^{1} \widehat{c}^{k}l_{i}l'_j \diff \xi\, u_{h,i}^{k} + \sum_{i=0}^{N}\int_{-1}^{1} l_{i}l'_j \diff \xi\, \widehat{f}_{\text{diff}\,h,i}^{k}. \label{eqn:integral4}
\end{align} 
Substituting the integrals (\ref{eqn:integral1}), (\ref{eqn:integral2}), (\ref{eqn:integral3}), and (\ref{eqn:integral4}), we get the following:

\begin{align}
\sum_{i=0}^{N}\left\{\frac{\Delta x_{k}}{2} \langle l_{i},l_{j} \rangle \frac{\diff u_{h,i}^{k}}{\diff t} + \left( \frac{\Delta x_{k}}{2\eta_1} \langle \chi^kl_{i},l_{j} \rangle - \langle \widehat{c}^{k}l_{i},l'_{j}\rangle \right) u_{h,i}^k - \langle l_{i},l'_j \rangle \widehat{f}_{\text{diff}\,h,i}^{k}\right\} = - \left.l_{j}\flux^k\right|_{-1}^{1},
\end{align}
for $k = 1,2,\ldots,K$ and $j = 0,1,\ldots,N$ where the inner product of the given functions $a = a(\xi)$ and $b = b(\xi)$ is defined as follows:

\begin{align}
\langle a,b\rangle \coloneqq \int_{-1}^{1}ab\,\diff\xi.
\end{align}
Additionally, the VP diffusive flux involves the derivative of $u$ and must be discretized consistently with the rest of the scheme. If we write Equation (\ref{eqn:VPdiffusiveflux}) in weak form, 

\begin{align}
\sum_{k=1}^{K}\left\{\frac{\Delta x_{k}}{2} \int_{-1}^{1}\left(\widehat{f}_{\text{diff}} + \widehat{\nu}\frac{2}{\Delta x_{k}}\frac{\partial u}{\partial \xi} \right) \psi\, \diff \xi \right\} = 0.
\end{align}
and repeat the interpolating and integration-by-part procedures, we get:

\begin{align}
\sum_{i=0}^{N} \left\{ \frac{\Delta x_{k}}{2} \langle l_{i},l_{j}\rangle \widehat{f}^{k}_{\text{diff}\,h,i} - \langle \widehat{\nu}^kl_{i},l'_{j} \rangle u^{k}_{h,i} \right\} = \left.-l_{j}\widehat{\nu}^k\numU^{k}\right|_{-1}^{1},
\end{align}
for $k = 1,2,\ldots,K$ and $j = 0,1,\ldots,N$ where $\numU^{k}$ is an analogy of the numerical flux for the solution, that is,

\begin{align}
& \numU_1^k \coloneqq \numU(u_{h,N}^k,u_{h,0}^{k+1}), \\
& \numU_{-1}^k \coloneqq \numU(u_{h,N}^{k-1},u_{h,0}^k).
\end{align}
The computation of the inner products is done via Gaussian quadrature:

\begin{align}
& \langle l_i,l_j\rangle \approx \sum_{m=0}^{N}w_ml_i(\xi_m)l_j(\xi_m) = w_j\delta_{ij}, \\
& \langle l_i,l'_j\rangle \approx \sum_{m=0}^{N}w_ml_i(\xi_m)l'_j(\xi_m) = w_il_j'(\xi_i), \\
& \langle \chi^kl_{i},l_{j} \rangle \approx \sum_{m=0}^{N}w_m\chi^k_ml_i(\xi_m)l_j(\xi_m) = w_j\chi^k_j\delta_{ij}, \\
& \langle \widehat{c}^{k}l_{i},l'_{j} \rangle \approx \sum_{m=0}^{N}w_m\widehat{c}^{k}_{m}l_i(\xi_m)l'_j(\xi_m) = w_i\widehat{c}^{k}_{i}l_j'(\xi_i), \\
& \langle \widehat{\nu}^kl_{i},l'_{j} \rangle \approx \sum_{m=0}^{N}w_m\widehat{\nu}^{k}_{m}l_i(\xi_m)l'_j(\xi_m) = w_i\widehat{\nu}^{k}_{i}l_j'(\xi_i),
\end{align}
where $w_{m}$ are the Gauss-Lobatto weights $(\sum_{m=0}^{N} w_m = 2)$ and $\chi^k_{m} = \chi^k(\xi_{m},t)$, $\widehat{c}^{k}_{m} = \widehat{c}^{k}(\xi_{m},t) = c + \chi^k_{m}/\eta_2$, and $\widehat{\nu}^{k}_{m} = \widehat{\nu}^{k}(\xi_{m},t) = \nu - \chi^k_{m}/\eta_3$. When all of them are combined, Equations (\ref{eqn:DGSEM}) are obtained.

Finally, the last stage of a DGSEM is the calculation of $\flux$ and $\numU$ to reproduce the physics of advection and diffusion. A variety of fluxes are available for DG, and most are summarized by Arnold et al. \cite{Arnold2002Elliptic}. Here, we use a unifying function: 

\begin{align}\label{eqn:Wfunction}
W^{k}_{\pm}(a,b;\lambda) \coloneqq \{\!\{ab\}\!\}^{k}_{\pm} - \frac{1}{2}\lambda[\![|a|b]\!]^{k}_{\pm},  
\end{align}
for $\lambda\in\mathbb{R}$. If $\lambda = 0$ the discretization becomes a central scheme; $\lambda = -1$, upwind; $\lambda = 1$, downwind. The subscript $``+"$ means the right boundary element and $``-"$ the left boundary element. We also denote $\{\!\{\cdot\}\!\}$ as the averaging operator:

\begin{align}
\{\!\{a\}\!\}^{k}_{+} \coloneqq \frac{a_{h,N}^{k} + a_{h,0}^{k+1}}{2},\quad \{\!\{a\}\!\}^{k}_{-} \coloneqq \frac{a_{h,0}^{k} + a_{h,N}^{k-1}}{2}, 
\end{align}
and $[\![\cdot]\!]$ as the jump operator:

\begin{align}
[\![a]\!]^{k}_{+} \coloneqq a_{h,N}^{k} - a_{h,0}^{k+1},\quad [\![a]\!]^{k}_{-} \coloneqq a_{h,N}^{k-1} - a_{h,0}^{k}. 
\end{align}
Once these operators are defined, we divide the numerical flux into an advective term and a diffusive term: $\flux = \flux_{\text{adv}} + \flux_{\text{diff}}$. The computation of the advective numerical flux is as follows:

\begin{align}
& \flux_{1,\text{adv}}^k = W^{k}_{+}(\widehat{c},u;\alpha), \\
& \flux_{-1,\text{adv}}^k = W^{k}_{-}(\widehat{c},u;\alpha),
\end{align}
and the diffusive numerical flux is: 

\begin{align}
& \flux_{1,\text{diff}}^k = W^{k}_{+}(1,\widehat{f}_{\text{diff}};\beta), \\
& \flux_{-1,\text{diff}}^k = W^{k}_{-}(1,\widehat{f}_{\text{diff}};\beta).
\end{align}
The values of $\widehat{f}_{\text{diff}}$ at the boundary elements are computed with:

\begin{align}
& \widehat{f}_{\text{diff}\,h,0}^{k} = W^{k}_{-}\left(\widehat{\nu},\frac{u}{\Delta x_k};\gamma\right), & & \widehat{f}_{\text{diff}\,h,N}^{k-1} = W^{k}_{+}\left(\widehat{\nu},\frac{u}{\Delta x_{k-1}};\gamma\right), \\
& \widehat{f}_{\text{diff}\,h,N}^{k} = W^{k}_{+}\left(\widehat{\nu},\frac{u}{\Delta x_k};\gamma\right), & & \widehat{f}_{\text{diff}\,h,0}^{k+1} = W^{k}_{-}\left(\widehat{\nu},\frac{u}{\Delta x_{k+1}};\gamma\right).
\end{align}
Finally, $\numU$ is computed as

\begin{align}
& \numU_{1}^k = W^{k}_{+}(1,u;\delta), \\ 
& \numU_{-1}^k = W^{k}_{-}(1,u;\delta). 
\end{align}
 BR1 is recovered by setting $\alpha = -1$ and $\beta = \gamma = \delta = 0$, while LDG is obtained by setting $\alpha = \gamma = -1$ and $\beta = -\delta = -1$. The weights $\mathfrak{f}$ and $\mathfrak{g}$ of (\ref{eqn:weightsfandg}) are obtained by finding the $u$s at $x_{k-1}$ and $x_{k}$ from the function $W$ described in Equation (\ref{eqn:Wfunction}).

\section{Non-trivial solutions}\label{appx:Nontrivial}

In all the cases, the considered element is inside the solid region. The first case is related to an inviscid problem without a second derivative penalty term or a viscid problem with $\eta_3 = 1/\nu$. The second case includes second derivatives and is therefore more general.

\subsection*{Case 1: Problem with $\widehat{\nu} = 0$} 
The parameters $TE$ and $HOT$ are listed in Tables \ref{tabappx:j_Reduced_PDE|Case1} and \ref{tabappx:NumericalSource|Case1}. The main findings include the following: 

\begin{align}
& \widetilde{c}^k_{j} = \Zhe_{j}^{(1)k} = \widehat{c},\quad \forall j = 0,1,2 \\
& \widetilde{\nu}^k_{j} = -\frac{1}{2}\Zhe_{j}^{(2)k} = 0,\quad \forall j = 0,1,2 \\
& \Zhe_{1}^{(2p)k} = 0, \quad p\in\mathbb{N}
\end{align}

In total, we have five unknowns ($\widehat{c},\,\mathfrak{f}^{k-1}_{2},\,\mathfrak{f}^{k}_{0},\,\mathfrak{f}^{k}_{2},\,\mathfrak{f}^{k+1}_{0}$) and, therefore, a determined system would be:

\begin{align}
\begin{cases}
\mathfrak{f}_{2}^{k-1}u_{h,2}^{k-1} + \left(\mathfrak{f}_{0}^{k} - \widehat{c}\right)u_{h,0}^{k} = 0, \\
\mathfrak{f}_{0}^{k+1}u_{h,0}^{k+1} + \left(\mathfrak{f}_{2}^{k} + \widehat{c}\right)u_{h,2}^{k} = 0, \\
\Zhe_{0}^{(3)k} = 0, \\
\Zhe_{1}^{(3)k} = 0, \\
\Zhe_{2}^{(3)k} = 0.
\end{cases}
\end{align}
whose errors are $TE_j^k\sim\Zhe_{j}^{(4)k}\sim\mathcal{O}\left(\Delta x_k^3\right)$ for $j = 0,1,2$ within $\Omega_\text{s}$. However, the unique solution of the system is the trivial one. If we leave $\eta_2$ free, the system that determines the numerical flux weights becomes 

\begin{align}
\begin{cases}
\mathfrak{f}_{2}^{k-1}u_{h,2}^{k-1} + \left(\mathfrak{f}_{0}^{k} - \widehat{c}\right)u_{h,0}^{k} = 0, \\
\mathfrak{f}_{0}^{k+1}u_{h,0}^{k+1} + \left(\mathfrak{f}_{2}^{k} + \widehat{c}\right)u_{h,2}^{k} = 0, 
\end{cases}
\end{align}
being $TE_j^k \sim \Zhe_{j}^{(3)k} \sim \mathcal{O}\left(\Delta x_k^2\right)$ for $j = 0,1,2$. A non-trivial solution would be:

\begin{align}
\mathfrak{f}_{2}^{k-1} = \mathfrak{f}_{0}^{k+1} = 0,\quad \mathfrak{f}_{0}^{k} = -\mathfrak{f}_{2}^{k} = \widehat{c},
\end{align} 
Alternatively, if the upwiding numerical flux is the solution, i.e. 

\begin{align}
\mathfrak{f}_{0}^{k} = \mathfrak{f}_{0}^{k+1} = 0,\quad\mathfrak{f}_{2}^{k-1} = -\mathfrak{f}_{2}^{k} = \widehat{c}. 
\end{align}
Then the system becomes:

\begin{align}
\widehat{c}\left(u_{h,2}^{k-1} - u_{h,0}^{k}\right) = 0, 
\end{align}
whose truncation error leads to:

\begin{align}
\begin{cases}
TE_j^k \sim \mathcal{O}\left(\Delta x_k^2\right),\enskip \forall j = 0,1,2, & \text{If } u_{h,2}^{k-1} = u_{h,0}^{k} \\
TE_0^k \sim \mathcal{O}\left(\Delta x_k^0\right),\enskip TE_1^k, TE_2^k \sim \mathcal{O}\left(\Delta x_k^2\right), & \text{If } u_{h,2}^{k-1} \neq u_{h,0}^{k}, 
\end{cases}
\end{align}
Setting $u_{h,2}^{k-1} = u_{h,0}^{k}$ is very similar to using a Continuous Galerkin (CG) method. A downwind numerical flux mimics the results of upwinding, but for the right-hand boundary element. Other numerical fluxes leave:

\begin{align}
TE_0^k,TE_2^k \sim \mathcal{O}\left(\Delta x_k^0\right),\quad TE_1^k \sim \mathcal{O}\left(\Delta x_k^2\right).
\end{align}

{
\renewcommand{\arraystretch}{2}
\begin{table}[!htpb]
\centering
% \begin{adjustbox}{angle=0}
\begin{tabular}{cccc}
\hline
$j$ & $\xi_j$ & $\widetilde{r}_{j}^{k}$ & $\Zhe_{j}^{(m)k}$        \\ \hline
$0$ & $-1$    & $\dfrac{6}{\Delta x_k}\widehat{c}$                     & $\dfrac{2^{2-m} - 1}{\Delta x_k^{1-m}}\widehat{c}$        \\ 
$1$ & $0$     & $0$                     & $- 2^{-m}\dfrac{(-1)^{m} - 1}{\Delta x_k^{1-m}}\widehat{c}$ \\ 
$2$ & $1$     & $-\dfrac{6}{\Delta x_k}\widehat{c}$                     & $(-1)^{m}\dfrac{1 - 2^{2-m}}{\Delta x_k^{1-m}}\widehat{c}$   \\ \hline
\end{tabular}
% \end{adjustbox}
\caption{The reaction parameter and the coefficient $\Zhe$ in the modified equations for a three-point GL grid and a problem with $\widehat{\nu} = 0$.}
\label{tabappx:j_Reduced_PDE|Case1}
\end{table}
}

{
\renewcommand{\arraystretch}{2}
\begin{table}[!htpb]
\centering
% \begin{adjustbox}{angle=0}
\begin{tabular}{cc}
\hline
$j$ & $\numS_{j}^{k}$ \\ \hline
$0$ & $\dfrac{6}{\Delta x_k}\left(\mathfrak{f}_{2}^{k-1}u_{h,2}^{k-1} + \mathfrak{f}_{0}^{k}u_{h,0}^{k}\right)$              \\
$1$ & $0$              \\
$2$ & $-\dfrac{6}{\Delta x_k}\left(\mathfrak{f}_{0}^{k+1}u_{h,0}^{k+1} + \mathfrak{f}_{2}^{k}u_{h,2}^{k}\right)$              \\ \hline
\end{tabular}
% \end{adjustbox}
\caption{Numerical source in the DG source, $s^k_{\text{DG},j} = \numS^k_j - \widetilde{r}^k_{j}u^{k}_{h,j}$, for a problem with $\widehat{\nu} = 0$.}
\label{tabappx:NumericalSource|Case1}
\end{table}
} 

\subsection*{Case 2: Problem with $\widehat{\nu} \neq 0$}
In this second case (with second derivatives) we consider $\eta_3$ free. The parameters of $TE$ and $HOT$ are listed in Tables \ref{tabappx:j_Reduced_PDE|Case2} and \ref{tabappx:NumericalSource|Case2}. Again, we conclude that 

\begin{align}
& \widetilde{c}^k_{0} = \Zhe_{0}^{(1)k} = \widehat{c} - 4\frac{3-\sfrac{3}{2}\mathfrak{g}^k_2}{\Delta x_k}\widehat{\nu}, \\
& \widetilde{c}^k_{1} = \Zhe_{1}^{(1)k} = \widehat{c} + 3\frac{\mathfrak{g}^k_2 - \mathfrak{g}^k_0}{\Delta x_k}\widehat{\nu}, \\
& \widetilde{c}^k_{2} = \Zhe_{2}^{(1)k} = \widehat{c} + 4\frac{3-\sfrac{3}{2}\mathfrak{g}^k_0}{\Delta x_k}\widehat{\nu}, 
\end{align}
and

\begin{align}
& \widetilde{\nu}^k_{0} = -\frac{1}{2}\Zhe_{0}^{(2)k} = (4-3\mathfrak{g}^k_2)\widehat{\nu}, \\
& \widetilde{\nu}^k_{1} = -\frac{1}{2}\Zhe_{1}^{(2)k} = -\frac{1}{4}\left(2+3\left(\mathfrak{g}^k_0+\mathfrak{g}^k_2\right)\right)\widehat{\nu}, \\
& \widetilde{\nu}^k_{2} = -\frac{1}{2}\Zhe_{2}^{(2)k} = (4-3\mathfrak{g}^k_0)\widehat{\nu}, 
\end{align}
Since $\widetilde{c}_j^k \neq \widehat{c}$ for $j=0,1,2$, the term $\widehat{c}^k_j - \Delta\xi\widetilde{c}^k_j$ in the truncation error should be suppressed since it is $\mathcal{O}(\Delta x_k^{-1})$. The choice is $\mathfrak{g}^k_2=\mathfrak{g}^k_0 = 2$. However, $\widehat{\nu}^k_j - \Delta\xi\widetilde{\nu}^k_j \neq 0$ for $j=0,1,2$ and therefore $TE_{j}^k \sim \mathcal{O}\left(\Delta x_k^0\right)$. If we want to find an optimal value of $\eta_3$ to increase the order of the scheme, we come to the conclusion that $\widehat{\nu} = 0$, but this case was already discussed previously. Keeping $\mathfrak{g}^k_2=\mathfrak{g}^k_0 = 2$ and $\eta_3$ free, $\mathfrak{g}^{k - 1}_2 = \mathfrak{g}^{k + 1}_0 = 0$ kills $s^k_{\text{DG},1}$, to have $s^k_{\text{DG},0}, s^k_{\text{DG},2} = 0$,

\begin{align}
\begin{cases}
\mathfrak{f}_{2}^{k-1}u_{h,2}^{k-1} + \left(\mathfrak{f}_{0}^{k} - \widehat{c} - \dfrac{4}{\Delta x_k}\widehat{\nu}\right)u_{h,0}^{k} = 0, \\
\mathfrak{f}_{0}^{k+1}u_{h,0}^{k+1} + \left(\mathfrak{f}_{2}^{k} + \widehat{c} + \dfrac{4}{\Delta x_k}\widehat{\nu}\right)u_{h,2}^{k} = 0. 
\end{cases}
\end{align}
In this case, a solution of the system is as follows:

\begin{align}
\mathfrak{f}_{2}^{k-1} = \mathfrak{f}_{0}^{k+1} = 0,\quad \mathfrak{f}_{0}^{k} = -\mathfrak{f}_{2}^{k} = \widehat{c} + \frac{4}{\Delta x_k}\widehat{\nu},    
\end{align}
Additionally, if upwind in such a way that

\begin{align}
\mathfrak{f}_{0}^{k} = \mathfrak{f}_{0}^{k+1} = 0,\quad \mathfrak{f}_{2}^{k-1} = -\mathfrak{f}_{2}^{k} = \widehat{c} + \frac{4}{\Delta x_k}\widehat{\nu},
\end{align}  
the second equation of the system is met, but the first one becomes:

\begin{align}
\left(\widehat{c} + \frac{4}{\Delta x_k}\widehat{\nu}\right)\left(u_{h,2}^{k-1} - u_{h,0}^{k}\right) = 0.
\end{align}
To eliminate this term, a relation of $\eta$s is obtained, $\widehat{c} + (4/\Delta x_k)\widehat{\nu} = 0$, since in a DG method $u_{h,2}^{k-1} \neq u_{h,0}^{k}$. Note that in a CG method, it is not necessary to fill this relation, since the solution is continuous between elements. In all cases described previously, $TE_{j}^{k} \sim \mathcal{O}(\Delta x_k^0) $ for $j=0,1,2$.

{
\renewcommand{\arraystretch}{2}
\begin{table}[!htpb]
\centering
% \begin{adjustbox}{angle=0}
% \resizebox{\textwidth}{!}{
\begin{tabular}{cccc}
\hline
$j$ & $\xi_j$ & $\widetilde{r}_{j}^{k}$ & $\Zhe_{j}^{(m)k}$        \\ \hline
$0$ & $-1$    & $\dfrac{6}{\Delta x_k}\widehat{c} - 4\dfrac{6}{\Delta x_k^2}\widehat{\nu}$                     & $\dfrac{2^{2-m} - 1}{\Delta x_k^{1-m}}\widehat{c} - 4\dfrac{2^{2-m} + 1-\sfrac{3}{2}\mathfrak{g}_{2}^{k}}{\Delta x_k^{2-m}}\widehat{\nu}$        \\ 
$1$ & $0$     & $2\dfrac{6}{\Delta x_k^2}\widehat{\nu}$                     & $- 2^{-m}\dfrac{(-1)^{m} - 1}{\Delta x_k^{1-m}}\widehat{c} + 2^{1-m}\dfrac{(-1)^{m}\left(1+3\mathfrak{g}_{0}^{k}\right) + 1+3\mathfrak{g}_{2}^{k}}{\Delta x_k^{2-m}}\widehat{\nu}$ \\ 
$2$ & $1$     & $-\dfrac{6}{\Delta x_k}\widehat{c} - 4\dfrac{6}{\Delta x_k^2}\widehat{\nu}$                     & $(-1)^{m}\left(\dfrac{1 - 2^{2-m}}{\Delta x_k^{1-m}}\widehat{c} - 4\dfrac{2^{2-m} + 1-\sfrac{3}{2}\mathfrak{g}_{0}^{k}}{\Delta x_k^{2-m}}\widehat{\nu}\right)$   \\ \hline
\end{tabular}
% }
% \end{adjustbox}
\caption{The reaction parameter and the coefficient $\Zhe$ in the modified equations for a three-point GL grid and a problem with $\widehat{\nu} \neq 0$.}
\label{tabappx:j_Reduced_PDE|Case2}
\end{table}
}

{
\renewcommand{\arraystretch}{2}
\begin{table}[!htpb]
\centering
% \begin{adjustbox}{angle=0}
% \resizebox{\textwidth}{!}{
\begin{tabular}{cc}
\hline
$j$ & $\numS_{j}^{k}$ \\ \hline
$0$ & $\dfrac{2}{\Delta x_k}\left[3\left(\mathfrak{f}_{2}^{k-1}u_{h,2}^{k-1} + \mathfrak{f}_{0}^{k}u_{h,0}^{k}\right) - \dfrac{3}{\Delta x_k} \left(3\mathfrak{g}_{2}^{k-1}u_{h,2}^{k-1} + \mathfrak{g}_{0}^{k+1}u_{h,0}^{k+1} + \left(3\mathfrak{g}_{0}^{k} + \mathfrak{g}_{2}^{k}\right)u_{h,0}^{k}\right)\widehat{\nu}\right]$              \\
$1$ & $\dfrac{6}{\Delta x_k^2}\left[\mathfrak{g}_2^{k-1}u_{h,2}^{k-1} + \mathfrak{g}_0^{k+1}u_{h,0}^{k+1} + \left(\mathfrak{g}_0^k + \mathfrak{g}_2^k\right)u_{h,1}^k\right]\widehat{\nu}$              \\
$2$ & $\dfrac{2}{\Delta x_k}\left[-3\left(\mathfrak{f}_{2}^{k}u_{h,2}^{k} + \mathfrak{f}_{0}^{k+1}u_{h,0}^{k+1}\right) - \dfrac{3}{\Delta x_k}\left(\mathfrak{g}_{2}^{k-1}u_{h,2}^{k-1} + 3\mathfrak{g}_{0}^{k+1}u_{h,0}^{k+1} + \left(3\mathfrak{g}_{2}^{k} + \mathfrak{g}_{0}^{k}\right)u_{h,2}^{k}\right)\widehat{\nu}\right]$              \\ \hline
\end{tabular}
% }
% \end{adjustbox}
\caption{Numerical source in the DG source, $s^k_{\text{DG},j} = \numS^k_j - \widetilde{r}^k_{j}u^{k}_{h,j}$, for a problem with $\widehat{\nu} \neq 0$.}
\label{tabappx:NumericalSource|Case2}
\end{table}
}

A summary of all the conditions derived can be found in Table \ref{tab:nontrivial}.

\section{Connection between Fourier series and Taylor series}\label{appx:FourierAndTaylor}

Assume that we analyze the errors of the advection equation,

\begin{align}
    \frac{\partial u}{\partial t} + c\frac{\partial u}{\partial x} = 0
\end{align}
with a positive non-zero constant, $c > 0$. For example and simplicity, the semidiscretization using a central finite difference is 

\begin{align}\label{eqn:semidiscreteexample}
    \frac{\diff u_j}{\diff t} + c\frac{u_{j+1} - u_{j-1}}{2\Delta x} = 0
\end{align}
with $\Delta x = x_{j+1} - x_{j} = x_{j} - x_{j-1}$ and $u_j = u(x_j,t)$. From the point of view of the Taylor series, the nodal solutions are expanded as

\begin{align}
    u_{j\pm 1} = u_{j} \pm \Delta x\left.\frac{\partial u}{\partial x}\right|_{j} + \frac{\Delta x^2}{2!}\left.\frac{\partial^2 u}{\partial x^2}\right|_{j} \pm \frac{\Delta x^3}{3!}\left.\frac{\partial^3 u}{\partial x^3}\right|_{j} + \dots 
\end{align}
and, therefore, the modified equation is 

\begin{align}\label{eqn:modifiedequationTaylor}
    \frac{\partial u_j}{\partial t} + c\left.\frac{\partial u}{\partial x}\right|_{j} = -c\frac{\Delta x^2}{6}\left.\frac{\partial^3 u}{\partial x^3}\right|_{j} + \dots \sim \mathbcal{O}(\Delta x^2).
\end{align}
From the point of view of the Fourier analysis, we start with a trial solution in the form:

\begin{align}\label{eqn:trialsolution}
    u(x_j,t) = u_\omega(t)\exp(\mathbf{i}\omega x_{j}),
\end{align}
where $\mathbf{i}$ is the imaginary unit and $\omega$ the wavenumber. The semi-discrete equation (\ref{eqn:semidiscreteexample}) now becomes 

\begin{align}
    \frac{\diff u_\omega}{\diff t} + \mathbf{i}c\frac{\sin(\omega\Delta x)}{\Delta x}u_\omega = 0.
\end{align}
The expansion of $\sin(\omega\Delta x)/\Delta x$ around $w = 0$ yields the following.

\begin{align}\label{eqn:aproxspectralequation}
    \frac{\diff u_\omega}{\diff t} + \mathbf{i}c\omega u_\omega - \mathbf{i}c\frac{\Delta x^2}{6}\omega^3 u_\omega + \dots = 0.
\end{align}
With the sinusoidal function (\ref{eqn:trialsolution}) we can relate the derivative of $m$ order to the wavenumber, that is,

\begin{align}\label{eqn:derivtrialsolution}
    \left.\frac{\partial^m u}{\partial x^m}\right|_{j} = (\mathbf{i}\omega)^m u_\omega\exp(\mathbf{i}\omega x_{j}) = (\mathbf{i}\omega)^m u_j .
\end{align}
Therefore, substituting Equation (\ref{eqn:derivtrialsolution}) into Equation (\ref{eqn:aproxspectralequation}) we obtain the modified equation:

\begin{align}
    \frac{\partial u_j}{\partial t} + c\left.\frac{\partial u}{\partial x}\right|_{j} = -c\frac{\Delta x^2}{6}\left.\frac{\partial^3 u}{\partial x^3}\right|_{j} + \dots \sim \mathbcal{O}(\Delta x^2)
\end{align}
which is the same as the modified equation (\ref{eqn:modifiedequationTaylor}) derived from the Taylor series. A similar and comprehensive understanding of the discrete errors by Fourier analysis can be found in \cite{Vichnevetsky1982fourier}. This connection is not coincidental, indeed, for the Taylor series we can obtain the trial solution and vice versa, e.g.,

\begin{align}
\nonumber    u_{j\pm 1} &= \sum_{m=0}^{\infty}(\pm 1)^m\frac{\Delta x^m}{m!}\left.\frac{\partial^m u}{\partial x^m}\right|_{j} \\
\nonumber               &= \sum_{m=0}^{\infty}(\pm 1)^m\frac{\Delta x^m}{m!}(\mathbf{i}\omega)^mu_j \\
\nonumber               &= u_j\sum_{m=0}^{\infty}\frac{(\pm\mathbf{i}\omega\Delta x)^m}{m!} \\
\nonumber               &= u_j\exp(\pm\mathbf{i}\omega\Delta x) \\ 
                        &= u_\omega\exp(\mathbf{i}\omega x_{j\pm 1})
\end{align}
Although a Taylor series approximates a function point-wise and a Fourier series approximates a function globally, both methods yield the same modified equation and later the same truncation error. This achievement can facilitate extension to multidimensional problems using the Fourier framework; see \cite{Voth2004Fourier}.

\bibliographystyle{spmpsci}      % mathematics and physical sciences
\bibliography{refs}   % name your BibTeX data base

% \begin{thebibliography}{}
% %
% % and use \bibitem to create references. Consult the Instructions
% % for authors for reference list style.
% % Format for Journal Reference
% Author, Article title, Journal, Volume, page numbers (year)
% Author, Book title, page numbers. Publisher, place (year)
% \end{thebibliography}

\end{document}